%% file: Group_draw_attractiveness_equal_treatment_handball.tex
\documentclass[12pt]{article}
\usepackage[utf8]{inputenc}
\usepackage[british]{babel}
\usepackage{cmap}
\usepackage{lmodern}

\usepackage{amssymb, amsmath, amsthm}
\usepackage[a4paper,top=25mm,bottom=25mm,left=25mm,right=25mm]{geometry}
\usepackage{ragged2e}

\usepackage{authblk} % for headings
\usepackage{pifont}
\usepackage{graphicx}
\usepackage[dvipsnames,svgnames,table]{xcolor}
\usepackage[figuresright]{rotating}
\usepackage{xtab} % tackle the long tables
\usepackage{longtable} % tackle the long tables
\usepackage{multirow}
\usepackage{footnote}
\usepackage[stable]{footmisc}
\usepackage{chngpage} % allows for temporary adjustment of side margins
\usepackage{pdflscape} % landscape environment
\usepackage[nottoc,notlot,notlof]{tocbibind} % includes everything in ToC

\usepackage{pgfplots}
\pgfplotsset{every tick label/.append style={font=\footnotesize}}
\pgfplotsset{compat=1.18}
\usepackage{setspace}

\usepackage{array}
\newcolumntype{K}[1]{>{\centering\arraybackslash$}p{#1}<{$}}

\makesavenoteenv{tabular}
\usepackage{tabularx}
\usepackage{booktabs}
\usepackage{threeparttable}
\usepackage[referable]{threeparttablex} % footnotes in tabu
\newcolumntype{R}{>{\raggedleft\arraybackslash}X}
\newcolumntype{L}{>{\raggedright\arraybackslash}X}
\newcolumntype{C}{>{\centering\arraybackslash}X}
\newcolumntype{A}{>{\columncolor{gray!25}}C}
\newcolumntype{a}{>{\columncolor{gray!25}}c}

\newlength{\tablen}

\usepackage{dcolumn} % alignment to decimal points
\newcolumntype{.}{D{.}{.}{-1}}

\usepackage{tikz}
\usetikzlibrary{arrows, calc, matrix, patterns, positioning, trees}
\usepackage[semicolon]{natbib} % numbers,
\usepackage[hyphens]{xurl}
\usepackage{microtype}
\usepackage[justification=centering]{caption} % [justification=centering]

\usepackage[labelformat=simple]{subcaption}

\DeclareCaptionLabelFormat{parenthesis}{(#2)}
\makeatletter
\renewcommand\p@subfigure{\arabic{figure}.}
\makeatother

\DeclareCaptionLabelFormat{parenthesis}{(#2)}
\makeatletter
\renewcommand\p@subtable{\arabic{table}.}
\makeatother

\usepackage{enumitem}

\setlist[itemize]{leftmargin=2.5\parindent}
\setlist[enumerate]{leftmargin=2.5\parindent}

\usepackage{hyperref} % [hidelinks]
\hypersetup{
  colorlinks   = true,    		% Colours links instead of ugly boxes
  urlcolor     = blue,    		% Colour for external hyperlinks
  linkcolor    = blue,    		% Colour of internal links
  citecolor    = ForestGreen	% Colour of citations
}

\newenvironment{customlegend}[1][]{%
	\begingroup
	\csname pgfplots@init@cleared@structures\endcsname
	\pgfplotsset{#1}%
    }{%
	\csname pgfplots@createlegend\endcsname
	\endgroup
    }%
\def\addlegendimage{\csname pgfplots@addlegendimage\endcsname}

\theoremstyle{plain}
\theoremstyle{definition}
\newtheorem{example}{Example}%[section]

\theoremstyle{remark}

\makeatletter
\let\@fnsymbol\@alph
\makeatother

\def\keywords{\vspace{.5em} % Add keywords
{\noindent \textit{Keywords}: }}

\def\AMS{\vspace{.5em} % Add keywords
{\noindent \textbf{\emph{MSC} class}: }}

\def\JEL{\vspace{.5em} % Add keywords
{\noindent \textbf{\emph{JEL} classification number}: }}

\title{Attractiveness and equal treatment in a group draw}

\author{\href{https://sites.google.com/view/laszlocsato}{L\'aszl\'o Csat\'o}\thanks{~Corresponding author. Email: \emph{laszlo.csato@sztaki.hun-ren.hu} \newline
Institute for Computer Science and Control (SZTAKI), Hungarian Research Network (HUN-REN), Laboratory on Engineering and Management Intelligence, Research Group of Operations Research and Decision Systems, Budapest, Hungary \newline
Corvinus University of Budapest (BCE), Institute of Operations and Decision Sciences, Department of Operations Research and Actuarial Sciences, Budapest, Hungary}
$\qquad \qquad \qquad \qquad$
\href{https://sites.google.com/view/doragretapetroczy}{D\'ora Gr\'eta Petr\'oczy}\thanks{~E-mail: \emph{apetroczy@metropolitan.hu} \newline
Budapest Metropolitan University, MNB Institute, Budapest, Hungary}}

\date{\today}

\def\Dedication{
{\noindent
``\emph{As a historian it is my duty to stress what we do not know.}''
}
\vspace{0.25cm}

\flushright
\noindent (Adrian Goldsworthy: \emph{Vindolanda})

\vspace{1cm} 
\justify }

\begin{document}

\maketitle

\thispagestyle{empty}
\Dedication

\begin{abstract}
\noindent
National teams from different continents can play against each other only in afew sports competitions. Therefore, a reasonable aim is maximising the number of intercontinental games in world cups, as done in basketball and football, in contrast to handball and volleyball. However, this objective requires additional draw constraints that imply the violation of equal treatment. In addition, the standard draw mechanism is non-uniformly distributed on the set of valid assignments, which may lead to further distortions.
Our paper analyses this novel trade-off between attractiveness and fairness through the example of the 2025 World Men's Handball Championship. We introduce a measure of inequality, which enables considering 32 sets of reasonable geographical restrictions to determine the Pareto frontier. The proposed methodology can be used by policy-makers to select the optimal set of draw constraints.


\AMS{62-08, 90-10, 90B90, 91B14}
% Computational methods for problems pertaining to statistics
% Mathematical modeling or simulation for problems pertaining to operations research and mathematical programming
% Case-oriented studies in operations research
% Social Choice

\JEL{C44, C63, Z20}
% Operations Research, Statistical Decision Theory
% Computational Techniques, Simulation Modeling 
% Sports Economics, General
\end{abstract}

\clearpage

\section{Introduction} \label{Sec1}

Attractiveness and fairness are among the most important criteria for evaluating the rules of a sports competition \citep{Csato2021a, DevriesereCsatoGoossens2024}. Tournament organisers face a difficult choice if these goals can only be improved at the expense of each other. Our paper analyses such a trade-off between the number of intra-continental matches and equal treatment in a group draw.

Geographic separation could be important for several reasons.
First, few sports competitions allow teams from different continents to play against each other. Hence, the uniqueness of a tournament is enhanced if the number of inter-continental games is maximised in the group stage \citep{Guyon2015a, LalienaLopez2019}.
Second, rating methods usually become more accurate and reliable if national teams from different continents play more matches against each other \citep{LasekGagolewski2018, SziklaiBiroCsato2022}.
Third, the qualification system is often distorted: from geographical zones where the particular sport is not so developed, weaker teams can qualify than from continents with several strong countries. This phenomenon is well-documented in football \citep{Csato2023c, KrumerMoreno-Ternero2023, CsatoKissSzadoczki2025} but is also relevant to other sports. For example, in the preliminary group stage of the 2025 IHF Men's World Championship, European and non-European teams played 26 matches but only three were won by non-European teams. Consequently, restricting the number of European teams in a group would be crucial to creating groups of roughly equal strength.

The balancedness of the groups is usually achieved by a seeding system: the teams are divided into pots according to their strength (an exogenous ranking based on historical performances), and each group contains at most one team from each pot. Therefore, equal treatment requires a team drawn from a given pot to have \emph{a priori} the same chance of playing against any team drawn from another pot. Otherwise, the draw rules may imply a higher probability of facing a long-term rival, which would be unfair. The importance of equal treatment is increased by the recent finding that nontransitivity can persist over a long time in sports \citep{vanOurs2024}.

Surprisingly, the governing bodies of major sports apply different policies regarding intra-continental matches in their tournaments. Both the group draw of recent FIBA World Cups \citep{FIBA2019, FIBA2023} and FIFA World Cups \citep{FIFA2017c, FIFA2022a}, the world championships in basketball and football, respectively, have imposed constraints on the set of feasible assignments: no group could have contained more than one team from any continent except for Europe, while lower and upper bounds have existed for the number of European teams.
On the other hand, analogous restrictions have not been used in the 2025 IHF Men's World Handball Championship \citep{IHF2025} and the 2025 FIVB Men's Volleyball World Championship \citep{Volleyball2024}.

This discrepancy has inspired the current research. To address a gap in the previous literature, we will explore a novel trade-off between the number of unattractive intra-continental matches and the extent of violating equal treatment by comparing different sets of draw constraints. The analysis will be conducted for the 2025 IHF World Men's Handball Championship draw. The complexity of the problem is increased by the standard procedure used to satisfy the draw constraints, which is non-uniformly distributed on the set of feasible assignments, potentially implying further inequality distortions.

Our main contributions can be summarised as follows:
\begin{itemize}
\item
First in the literature, a trade-off between attractiveness and equal treatment is revealed and studied in the draw of a sports tournament;
\item
A measure based on the Herfindahl--Hirschman index is proposed to quantify the degree to which the principle of equal treatment is violated (Section~\ref{Sec33});
\item
By examining 32 sets of reasonable draw constraints, the standard mechanism for group draw is shown to struggle with a restriction that allows more than one team from the same geographic zone to play in a group (Figure~\ref{Fig2});
\item
Eight sets of draw constraints are uncovered to be suboptimal as they imply a dominated combination in both dimensions, attractiveness and fairness (Figure~\ref{Fig3}).
\end{itemize}

The paper is structured as follows.
Section~\ref{Sec2} gives an overview of the existing literature. The case study from handball and our methodology are detailed in Section~\ref{Sec3}. Section~\ref{Sec4} presents and discusses the results, while Section~\ref{Sec5} provides concluding remarks.

\section{Related literature} \label{Sec2}

Analysing trade-offs is a standard research approach in the field of tournament design.
\citet{GoossensYiVanBulck2020} discuss three fairness issues (consecutive home games, carry-over effect, number of rest days between consecutive games) in sports timetabling.
\citet{LasekGagolewski2018} and \citet{SziklaiBiroCsato2022} investigate the relationship between the number of matches played and the accuracy of the ranking.
\citet{Csato2022a} reveals a trade-off between draw restrictions and the probability of a game that creates perverse incentives to lose for a team.
\citet{SauerCsehLenzner2024} propose alternative pairing mechanisms for Swiss-system tournaments to optimise various objective functions.
The recent survey \citet{DevriesereCsatoGoossens2024} presents several further examples.

The problem of group draw has also received serious attention by the academic community.
\citet{Guyon2015a} outlines the main criteria for the FIFA World Cup draw (balance, even/uniform distribution, fairness, geographic separation, randomness, tractability) and identifies some flaws in the 2014 FIFA World Cup draw.
A preliminary version of this paper, \citet{Guyon2014a} suggests three alternative draw procedures, one of which inspired FIFA to introduce the so-called \emph{Skip mechanism} (see Section~\ref{Sec32}) for the 2018 FIFA World Cup draw \citep{Guyon2018d}.

One line of research deals with balancing the strengths of the groups.
\citet{LalienaLopez2019} consider draw systems for tournaments with 32 teams and eight groups under seeding and geographical restrictions. Their first method creates perfectly balanced groups by listing all feasible solutions, however, its practical implementation remains challenging. Their second method is a heuristic, which does not guarantee perfect balance but provides good results for the 2014 FIFA World Cup draw.
\citet{CeaDuranGuajardoSureSiebertZamorano2020} construct a mixed integer linear programming model to generate group assignments satisfying geographical separation and ensuring balance between the strengths of the groups.
\citet{Csato2023d} demonstrates that the 2022 FIFA World Cup draw has not balanced the groups to the extent possible because the winners of the play-offs were assigned to the weakest pot.
\citet{LalienaLopez2025} start from the ideas of \citet{LalienaLopez2019} to develop a more efficient algorithm for groups of three teams where the strength of a group is constrained by exogenously given lower and upper bounds.

A further research area is analysing the bias of real-world draw mechanisms compared to a uniform draw.
Between the 2003/04 and the 2023/24 seasons, the UEFA Champions League Round of 16 draw matched eight group winners and eight runners-up such that teams from the same national association and the same group could not played against each other. The draw was implemented by the so-called \emph{Drop mechanism}. \citet{Kiesl2013} computes the distortion in the 2012/13 season and verifies the existence of a feasible assignment by Hall's marriage theorem \citep{Hall1935}. Both issues are discussed by \citet{WallaceHaigh2013}, too. \citet{KlossnerBecker2013} prove that the Drop mechanism cannot be uniformly distributed in this setting and present the quite substantial financial consequences of small probability differences. The algebraic structure of the permutation describing the outcome of the draw is found to have a powerful impact on its probability. \citet{GuyonMeunier2023} explain the calculation of exact draw probabilities by using well-known maximum matching algorithms.
Last but not least, the Drop mechanism seems to be close to a constrained-best in the UEFA Champions League Round of 16 draw according to the results of \citet{BoczonWilson2023}.

UEFA has introduced an incomplete round-robin format for its club competitions from the 2024/25 season \citep{Gyimesi2024}. The draw procedure of the new league phase of the UEFA Champions League is examined by \citet{GuyonMeunierBenSalemBuchholtzerTanre2024}.

The distortion of the Skip mechanism has also been studied.
\citet{RobertsRosenthal2024} discuss how it can be simulated and quantifies the bias of the 2022 FIFA World Cup draw. In the 2018 FIFA World Cup draw, \citet{Csato2025c} uncovers the effect of the draw order on the level of unfairness. \citet{Csato2024i} explains the conclusions from these numerical calculations in a theoretical framework. \citet{Csato2024h} compares the performance of the Drop and Skip mechanisms for bipartite graphs, including those corresponding to the historical seasons of the UEFA Champions League Round of 16.

The novelty of the current paper resides in our focus on two requirements for group draws simultaneously: fairness (equal treatment) and geographic separation. We do not know of any study that investigates multiple goals for draw procedures. In addition, previous papers have always considered a given set of draw constraints, but they are allowed to vary in the following according to the preferences of the organiser.
%In addition, uniform distribution will be considered by comparing the Skip mechanism with a uniform draw.

\section{Methodology} \label{Sec3}

The rules of the 2025 IHF World Men's Handball Championship draw and some reasonable geographical constraints are described in Section~\ref{Sec31}. Section~\ref{Sec32} presents draw procedures that can be used if further restrictions are imposed. Section~\ref{Sec33} defines our measure to quantify the extent to which the principle of equal treatment is violated.

\subsection{The 2025 IHF World Men's Handball Championship draw} \label{Sec31}

\begin{table}[t!]
  \centering
  \caption{The geographical composition of the seeding pots \\ in the 2025 IHF Men's World Handball Championship draw}
  \label{Table1}
  \rowcolors{3}{}{gray!20}
    \begin{tabularx}{0.8\textwidth}{l CCCCC} \toprule
    Confederation & Pot 1 & Pot 2 & Pot 3 & Pot 4 & Sum \\ \bottomrule
    Africa & 1     & 0     & 1     & 3     & 5 \\
    Asia  & 0     & 0     & 2     & 2     & 4 \\
    Europe & 7     & 8     & 2     & 1     & 18 \\
    North America & 0     & 0     & 1     & 1     & 2 \\
    South America & 0     & 0     & 2     & 1     & 3 \\ \toprule
    \end{tabularx}
\end{table}

The 2025 IHF World Men's Handball Championship has been contested by 32 teams. In particular, five teams have arrived from Africa, four from Asia, 18 from Europe, two from the North American and Caribbean Confederation (in the following, North America), and three from the South and Central American Confederation (in the following, South America).
Prior to the draw, the teams have been seeded into four pots, whose geographical composition is shown in Table~\ref{Table1}. Each group has contained one team from each pot.

\begin{table}[t!]
  \centering
  \caption{Group composition in the preliminary round of \\ the 2025 IHF Men's World Handball Championship}
  \label{Table2}
  \rowcolors{3}{gray!20}{}
    \begin{tabularx}{1\textwidth}{l CCCC CCCC C} \toprule \hiderowcolors
    \multirow{2}{*}{Confederation} & \multicolumn{8}{c}{Group} & \multirow{2}{*}{Sum} \\    	
          & A & B & C & D & E & F & G & H & \\ \bottomrule \showrowcolors
    Africa & 0     & 2     & 0     & 1     & 0     & 0     & 1     & 1     & 5 \\
    Asia  & 0     & 0     & 2     & 0     & 0     & 1     & 0     & 0     & 3 \\
    Europe & 4     & 2     & 2     & 3     & 2     & 2     & 2     & 1     & 18 \\
    North America & 0     & 0     & 0     & 0     & 1     & 0     & 1     & 1     & 3 \\
    South America & 0     & 0     & 0     & 0     & 1     & 1     & 0     & 1     & 3 \\ \toprule
    Unattractive games & 6     & 2     & 2     & 3     & 1     & 1     & 1     & 0     & 16 \\ \bottomrule
    \end{tabularx}
\end{table}

Groups A and B have played in Denmark, Groups E and F in Norway, and Groups C, D, G, H in Croatia. Therefore, Denmark (Pot 1) has automatically been assigned to Group B, Norway (Pot 1) to Group E, and Croatia (Pot 2) to Group H \citep{IHF2025}. Furthermore, the hosts have had the right to choose one team for the other groups that played in their country. The Danish Handball Association has placed Germany (Pot 1) in Group A, while the Norwegian Handball Federation has placed Sweden in Group F. Due to the decision of Croatia, Austria (Pot 2) has been assigned to Group C, Hungary (Pot 1) to Group D, and Slovenia (Pot 1) to Group G.
No further draw constraints have been applied. Table~\ref{Table2} reports the final group composition.

As we have argued in the Introduction, intra-confederation games are undesirable and therefore called \emph{unattractive} in the following. Their number could be reduced by draw restrictions analogous to the ones used in the FIBA Basketball World Cup and the FIFA World Cup.
In particular, the following constraints are considered:
\begin{itemize}
\item
Constraint A: two African teams cannot play in the same group;
\item
Constraint B: two Asian teams cannot play in the same group;
\item
Constraint C: two North American teams cannot play in the same group;
\item
Constraint D: two South American teams cannot play in the same group;
\item
Constraint E: each group should contain at least two and at most three European teams.
\end{itemize}
The aim of Constraints A--D is straightforward. Regarding Constraint E, since there are 18 European teams for eight groups, the average number of European teams in a group is between two and three. Thus, Constraint E minimises the number of matches between two European teams because 1/2/3/4 teams from the same continent in a group results in 0/1/3/6 unattractive matches, respectively.

\begin{table}[t!]
  \centering
  \caption{Sets of possible draw constraints in the \\ 2025 IHF Men's World Handball Championship draw}
  \label{Table3}
  \rowcolors{3}{}{gray!20}
\centerline{
    \begin{tabularx}{1.15\textwidth}{l CCccC} \toprule \hiderowcolors
    Continent & Africa & Asia  & North America & South America & Europe \\
    Teams in a group & At most 1 & At most 1 & At most 1 & At most 1 & 1 or 2 \\
    Scenario & Constraint A & Constraint B & Constraint C & Constraint D & Constraint E \\ \bottomrule \showrowcolors
    0     & \textcolor{BrickRed}{\ding{55}} & \textcolor{BrickRed}{\ding{55}} & \textcolor{BrickRed}{\ding{55}} & \textcolor{BrickRed}{\ding{55}} & \textcolor{BrickRed}{\ding{55}} \\
    1     & \textcolor{BrickRed}{\ding{55}} & \textcolor{BrickRed}{\ding{55}} & \textcolor{BrickRed}{\ding{55}} & \textcolor{BrickRed}{\ding{55}} & \textcolor{PineGreen}{\ding{52}} \\
    2     & \textcolor{BrickRed}{\ding{55}} & \textcolor{BrickRed}{\ding{55}} & \textcolor{BrickRed}{\ding{55}} & \textcolor{PineGreen}{\ding{52}} & \textcolor{BrickRed}{\ding{55}} \\
    3     & \textcolor{BrickRed}{\ding{55}} & \textcolor{BrickRed}{\ding{55}} & \textcolor{BrickRed}{\ding{55}} & \textcolor{PineGreen}{\ding{52}} & \textcolor{PineGreen}{\ding{52}} \\
    4     & \textcolor{BrickRed}{\ding{55}} & \textcolor{BrickRed}{\ding{55}} & \textcolor{PineGreen}{\ding{52}} & \textcolor{BrickRed}{\ding{55}} & \textcolor{BrickRed}{\ding{55}} \\
    5     & \textcolor{BrickRed}{\ding{55}} & \textcolor{BrickRed}{\ding{55}} & \textcolor{PineGreen}{\ding{52}} & \textcolor{BrickRed}{\ding{55}} & \textcolor{PineGreen}{\ding{52}} \\
    6     & \textcolor{BrickRed}{\ding{55}} & \textcolor{BrickRed}{\ding{55}} & \textcolor{PineGreen}{\ding{52}} & \textcolor{PineGreen}{\ding{52}} & \textcolor{BrickRed}{\ding{55}} \\
    7     & \textcolor{BrickRed}{\ding{55}} & \textcolor{BrickRed}{\ding{55}} & \textcolor{PineGreen}{\ding{52}} & \textcolor{PineGreen}{\ding{52}} & \textcolor{PineGreen}{\ding{52}} \\
    8     & \textcolor{BrickRed}{\ding{55}} & \textcolor{PineGreen}{\ding{52}} & \textcolor{BrickRed}{\ding{55}} & \textcolor{BrickRed}{\ding{55}} & \textcolor{BrickRed}{\ding{55}} \\
    9     & \textcolor{BrickRed}{\ding{55}} & \textcolor{PineGreen}{\ding{52}} & \textcolor{BrickRed}{\ding{55}} & \textcolor{BrickRed}{\ding{55}} & \textcolor{PineGreen}{\ding{52}} \\
    10    & \textcolor{BrickRed}{\ding{55}} & \textcolor{PineGreen}{\ding{52}} & \textcolor{BrickRed}{\ding{55}} & \textcolor{PineGreen}{\ding{52}} & \textcolor{BrickRed}{\ding{55}} \\
    11    & \textcolor{BrickRed}{\ding{55}} & \textcolor{PineGreen}{\ding{52}} & \textcolor{BrickRed}{\ding{55}} & \textcolor{PineGreen}{\ding{52}} & \textcolor{PineGreen}{\ding{52}} \\
    12    & \textcolor{BrickRed}{\ding{55}} & \textcolor{PineGreen}{\ding{52}} & \textcolor{PineGreen}{\ding{52}} & \textcolor{BrickRed}{\ding{55}} & \textcolor{BrickRed}{\ding{55}} \\
    13    & \textcolor{BrickRed}{\ding{55}} & \textcolor{PineGreen}{\ding{52}} & \textcolor{PineGreen}{\ding{52}} & \textcolor{BrickRed}{\ding{55}} & \textcolor{PineGreen}{\ding{52}} \\
    14    & \textcolor{BrickRed}{\ding{55}} & \textcolor{PineGreen}{\ding{52}} & \textcolor{PineGreen}{\ding{52}} & \textcolor{PineGreen}{\ding{52}} & \textcolor{BrickRed}{\ding{55}} \\
    15    & \textcolor{BrickRed}{\ding{55}} & \textcolor{PineGreen}{\ding{52}} & \textcolor{PineGreen}{\ding{52}} & \textcolor{PineGreen}{\ding{52}} & \textcolor{PineGreen}{\ding{52}} \\
    16    & \textcolor{PineGreen}{\ding{52}} & \textcolor{BrickRed}{\ding{55}} & \textcolor{BrickRed}{\ding{55}} & \textcolor{BrickRed}{\ding{55}} & \textcolor{BrickRed}{\ding{55}} \\
    17    & \textcolor{PineGreen}{\ding{52}} & \textcolor{BrickRed}{\ding{55}} & \textcolor{BrickRed}{\ding{55}} & \textcolor{BrickRed}{\ding{55}} & \textcolor{PineGreen}{\ding{52}} \\
    18    & \textcolor{PineGreen}{\ding{52}} & \textcolor{BrickRed}{\ding{55}} & \textcolor{BrickRed}{\ding{55}} & \textcolor{PineGreen}{\ding{52}} & \textcolor{BrickRed}{\ding{55}} \\
    19    & \textcolor{PineGreen}{\ding{52}} & \textcolor{BrickRed}{\ding{55}} & \textcolor{BrickRed}{\ding{55}} & \textcolor{PineGreen}{\ding{52}} & \textcolor{PineGreen}{\ding{52}} \\
    20    & \textcolor{PineGreen}{\ding{52}} & \textcolor{BrickRed}{\ding{55}} & \textcolor{PineGreen}{\ding{52}} & \textcolor{BrickRed}{\ding{55}} & \textcolor{BrickRed}{\ding{55}} \\
    21    & \textcolor{PineGreen}{\ding{52}} & \textcolor{BrickRed}{\ding{55}} & \textcolor{PineGreen}{\ding{52}} & \textcolor{BrickRed}{\ding{55}} & \textcolor{PineGreen}{\ding{52}} \\
    22    & \textcolor{PineGreen}{\ding{52}} & \textcolor{BrickRed}{\ding{55}} & \textcolor{PineGreen}{\ding{52}} & \textcolor{PineGreen}{\ding{52}} & \textcolor{BrickRed}{\ding{55}} \\
    23    & \textcolor{PineGreen}{\ding{52}} & \textcolor{BrickRed}{\ding{55}} & \textcolor{PineGreen}{\ding{52}} & \textcolor{PineGreen}{\ding{52}} & \textcolor{PineGreen}{\ding{52}} \\
    24    & \textcolor{PineGreen}{\ding{52}} & \textcolor{PineGreen}{\ding{52}} & \textcolor{BrickRed}{\ding{55}} & \textcolor{BrickRed}{\ding{55}} & \textcolor{BrickRed}{\ding{55}} \\
    25    & \textcolor{PineGreen}{\ding{52}} & \textcolor{PineGreen}{\ding{52}} & \textcolor{BrickRed}{\ding{55}} & \textcolor{BrickRed}{\ding{55}} & \textcolor{PineGreen}{\ding{52}} \\
    26    & \textcolor{PineGreen}{\ding{52}} & \textcolor{PineGreen}{\ding{52}} & \textcolor{BrickRed}{\ding{55}} & \textcolor{PineGreen}{\ding{52}} & \textcolor{BrickRed}{\ding{55}} \\
    27    & \textcolor{PineGreen}{\ding{52}} & \textcolor{PineGreen}{\ding{52}} & \textcolor{BrickRed}{\ding{55}} & \textcolor{PineGreen}{\ding{52}} & \textcolor{PineGreen}{\ding{52}} \\
    28    & \textcolor{PineGreen}{\ding{52}} & \textcolor{PineGreen}{\ding{52}} & \textcolor{PineGreen}{\ding{52}} & \textcolor{BrickRed}{\ding{55}} & \textcolor{BrickRed}{\ding{55}} \\
    29    & \textcolor{PineGreen}{\ding{52}} & \textcolor{PineGreen}{\ding{52}} & \textcolor{PineGreen}{\ding{52}} & \textcolor{BrickRed}{\ding{55}} & \textcolor{PineGreen}{\ding{52}} \\
    30    & \textcolor{PineGreen}{\ding{52}} & \textcolor{PineGreen}{\ding{52}} & \textcolor{PineGreen}{\ding{52}} & \textcolor{PineGreen}{\ding{52}} & \textcolor{BrickRed}{\ding{55}} \\
    31    & \textcolor{PineGreen}{\ding{52}} & \textcolor{PineGreen}{\ding{52}} & \textcolor{PineGreen}{\ding{52}} & \textcolor{PineGreen}{\ding{52}} & \textcolor{PineGreen}{\ding{52}} \\ \bottomrule
    \end{tabularx}
}
\end{table}

The five constraints imply 32 different scenarios, as each can be either imposed or ignored. In the following, they are denoted by integers from 0 to 31, see Table~\ref{Table3}.

\subsection{The implementation of a draw with restrictions} \label{Sec32}

If draw constraints exist, finding a valid assignment of the teams into groups is non-trivial. A fair draw procedure should be uniformly distributed, that is, each valid assignment should have the same probability of occurring \citep{Guyon2015a, KlossnerBecker2013}. This will be called the \emph{Uniform mechanism} in the following.
Its simplest implementation is provided by a rejection sampler \citep[Section~2.1]{RobertsRosenthal2024}:
\begin{enumerate}
\item
Select uniformly at random from the set of $6! \times 2! \times \left( 8! \right)^2 \approx 2.34 \times 10^{12}$ possible assignments, which is obvious;
\item
If \emph{all} draw constraints are satisfied, the chosen assignment is accepted;
\item
If at least one draw constraint is violated, the chosen assignment is rejected, and the procedure returns to Step 1.
\end{enumerate}

According to Section~\ref{Sec31}, the number of valid assignments is already reduced by the selection of the hosts since two teams from Pot 2 (Austria and Croatia) cannot play against six teams from Pot 1 (Denmark, Germany, Hungary, Norway, Slovenia, Sweden). In other words, Austria and Croatia should play against Egypt and France from Pot 1.
This restriction, imposed by the organiser, will not be examined in our study.

However, no group draw uses a rejection sampler in practice because
(1) it would be boring and unsuitable for a streamed TV show; and
(2) it would threaten transparency as it is impossible to check for the public that the draw procedure is not manipulated for the sake of certain teams \citep{Tijms2015}.
Usually, the \emph{Skip mechanism} is applied if the number of groups exceeds two, such as in the FIBA Basketball World Cup and the FIFA World Cup \citep{Csato2024h}.

The Skip mechanism is a sequential algorithm that works as follows:
\begin{itemize}
\item
The draw starts with Pot 1 and continues with Pot 2 until Pot 4.
\item
Each pot is emptied before proceeding to the next pot.
\item
For each pot, the teams drawn are allocated in ascending alphabetical order from Group A to Group H.
\item
If a draw condition applies or is anticipated to apply, the team drawn is allocated to the next available group in alphabetical order as indicated by the computer.
\end{itemize}
The Skip mechanism always outputs a feasible allocation if it exists. However, because it is non-uniform, some valid assignments are more likely to occur \citep{Csato2025c, RobertsRosenthal2024}.

This draw procedure is not as simple as it might appear at first glance. In particular, the number of options available to a team depends not only on its own attributes and those of the teams already drawn, but also on the attributes of the teams still to be drawn. Otherwise, a deadlock situation might arise when the teams still to be drawn cannot be assigned to the remaining empty slots with satisfying all constraints.
Let us see an illustration.

\begin{example} \label{Examp1}
Assume that six teams 1--6 are drawn from Pots I (teams 1--3) and II (teams 4--6). The restrictions exclude two pairs of teams to play against each other: teams 1 and 4, teams 3 and 6. The teams are drawn according to their numerical order. Thus, when the draw reaches Pot II, Group A contains team 1, Group B contains team 2, and Group C contains team 3.
Team 4 cannot be assigned to Group A because the pair (1,4) is prohibited, it is placed in Group B. Since no constraint applies for team 5, its ``natural'' place would be Group A. However, team 6 cannot be assigned to Group C because the pair (3,6) is also prohibited. Consequently, team 2 should be assigned to Group C even though Group A still contains an empty slot for a team drawn from Pot II.
\end{example}

Although the definition of the Skip mechanism is simple, it is surprisingly challenging to simulate with a computer program \citep[p.~563]{RobertsRosenthal2024}. No organiser provides a method for this purpose, however, appropriate backtracking algorithms are presented in \citet{Csato2025c} and \citet{RobertsRosenthal2024}.
A simulator of the Skip mechanism for the 2018 and 2022 FIFA World Cups is available at \url{http://probability.ca/fdraw/}.

The choice of the host nations, discussed in Section~\ref{Sec31}, can be easily added analogous to Constraints A--E: two countries from the set of the eight particular nations (the three hosts and the five teams chosen by them) cannot play in the same group. Then the number of possible scenarios for the Skip mechanism---if group labels are ignored---is $\left( 8! \right)^3 \approx 6.55 \times 10^{13}$ since the teams can be drawn in an arbitrary order from each pot.

The outcome of the Skip mechanism may depend on the order of the seeding pots \citep{Csato2025c}. This issue is not examined here, we always use the procedure with the numerical order from Pot 1 to Pot 4, which is likely optimal with respect to the distortions compared to the Uniform mechanism \citep{Csato2024i}.

Obviously, complete enumeration is excluded for both the Uniform and the Skip mechanisms due to the huge number of assignments: the probabilities are determined based on one million random feasible assignments. For the Uniform mechanism, this requires generating much more unconstrained assignments that are checked by a rejection sampler. The running time is several hours for each set of draw constraints and for each draw mechanism on a standard laptop.

\subsection{The inequality of a group draw} \label{Sec33}

The tournament organiser has two reasonable objective functions: maximising the number of inter-continental (which is equivalent to minimising the number of unattractive games) and minimising the inequality imposed by the draw constraints. The former can be directly quantified. The level of inequality is measured as follows.

Denote the number of pots by $m$ and the number of teams in a pot by $n$. The $n \times n$ \emph{doubly stochastic matrix} $\mathbf{P}^{(k, \ell)} = \left[ p_{ij}^{(k, \ell)} \right]$ contains the probability $p_{ij}$ of assigning two teams $i$ and $j$ drawn from two different Pots $k$ and $\ell$ to the same group. The average Herfindahl--Hirschman index \citep{Herfindahl1950, Hirschman1945}, a well-known measure of market competition and competitive balance in sports \citep{OwenOwen2022, OwenRyanWeatherston2007}, for each row and each column of the $m(m-1)/2$ doubly stochastic assignment matrices is:
\[
\hat{I} = \frac{2}{m(m-1)} \sum_{1 \leq k, \ell \leq m} \frac{1}{2n} \sum_{1 \leq i \leq n} \sum_{1 \leq j \leq n} \left( p_{ij}^{(k, \ell)} \right)^2.
\]

The minimum of $\hat{I}$, $1/n$, is reached if all entries of all assignment matrices equal $1/n$. In this case, each team faces the highest uncertainty in the draw since they will play against each possible opponent with an equal probability of $1/n$.

The maximum of $\hat{I}$, 1, is reached if all doubly stochastic assignment matrices contain only zeros and ones. In this case, the draw is fully deterministic and unexciting.

$\hat{I}$ is transformed to the unit interval $\left[ 0, 1 \right]$ to get a \emph{measure of inequality}:
\[
I = \frac{\hat{I} - 1/n}{1- 1/n}.
\]

\begin{example} \label{Examp2}
Take the official unconstrained draw (Scenario 0) when Constraints A--E are not imposed in the 2025 IHF Men's Handball Championship. All entries of matrices $\mathbf{P}^{(1,3)}$, $\mathbf{P}^{(1,4)}$, $\mathbf{P}^{(2,3)}$, $\mathbf{P}^{(2,4)}$, $\mathbf{P}^{(3,4)}$ can be verified to equal $1/8$. However, this does not hold for the assignment matrix $\mathbf{P}^{(1,2)}$. Here, six rows and six columns contain $1/6$ six times, while two rows and two columns (corresponding to Egypt, France, Austria, Croatia) contain $1/2$ twice. Consequently,
\[
\frac{1}{2 \cdot 8} \sum_{1 \leq i \leq 8} \sum_{1 \leq j \leq 8} \left( p_{ij}^{(1,2)} \right)^2 = \frac{1}{16} \cdot \left[ 12 \cdot 6 \cdot (1/6)^2 + 4 \cdot 2 \cdot (1/2)^2 \right] = \frac{1}{16} \cdot \left( 2 + 2 \right) = \frac{1}{4}.
\]
This implies $\hat{I} = 1/6 \cdot \left( 5 \cdot 1/8 + 1/4 \right) = 7/48$ and $I = \left( 7/48 - 6/48 \right) \cdot 8/7 = 1/42$.
\end{example}

According to our knowledge, only \citet{BoczonWilson2023} have quantified previously the violation of equal treatment in a group draw. Their measure is the average absolute difference in the probabilities across all pairwise comparisons that are not directly excluded by the constraints. However, this idea cannot be applied here because the set of games with a non-zero probability depends on the set of draw constraints imposed, which varies in our case study.

\section{Results} \label{Sec4}

\input{Figure1_unattractive_match_distribution}

Figure~\ref{Fig1} shows the distribution of intra-continental matches in the absence of draw constraints (Scenario 0). European teams play at least 12 games against each other, which is the minimum. The probability exceeds 25\% for 14 and 15 unattractive matches, as well as for at least 16. The chance that the number of these games is above 18 remains below 1\%. 
According to Table~\ref{Table2}, 16 intra-continental games have been played in the 2025 IHF Men's World Handball Championship, which could have been guaranteed to decrease by 25\% with Constraints A--E.

\input{Figure2_constraints_impact_feasibility}

Imposing draw constraints makes some assignments invalid as presented in Figure~\ref{Fig2}.
Constraint E, which appears in the set of constraints denoted by an odd number, is clearly the most restrictive. Adding Constraint E excludes almost 70\% of possible allocations, while Constraint C is violated only in one case out of eight. Naturally, increasing the number of constraints decreases the probability of a valid assignment. Constraints A--D imply that only one out of eight allocations is feasible, however, this proportion remains about 5\% even if all of them are required. In contrast, the probability of a valid draw was merely 1/161 and 1/560 in the 2018 \citep{Csato2025c} and the 2022 \citep{RobertsRosenthal2024} FIFA World Cups, respectively. Thus, the low frequency of a feasible solution cannot be an argument against using geographical restrictions in the 2025 Handball World Championship.

Figure~\ref{Fig2} also uncovers the extent to which the Skip mechanism increases the level of inequality compared to the fair Uniform mechanism. The main challenge is posed by Constraint E, which allows up to three European teams in a group. While Constraints A--D can be guaranteed by some prohibited team pairs, this is not true if more than one team from a given set (in our case, from Europe) can play in a group. Therefore, if Constraint E is required, inequality is increased by at least 20\% (30\%) if Constraint A is (not) imposed due to the imperfect draw procedure.

\input{Figure3_attractiveness_inequality_tradeoff}

Figure~\ref{Fig3} conveys the main message of our paper by showing 64 potential group draws---implied by 32 different sets of constraints and two draw mechanisms---in the space determined by the two objective functions. The average number of unattractive matches can vary between 12 and 14.62, while the inequality lies between $1/42 \approx 0.024$ (see Example~\ref{Examp2}) and $0.05$.
The trade-off can be summarised as follows:
\begin{itemize}
\item
The most efficient curve is given by the eight combinations of Constraints B--D that affect only the teams drawn from Pots 3 and 4. Along this line, the expected number of intra-continental matches can be reduced by almost one with increasing inequality by about 16\%. The Skip mechanism is not worse than the Uniform.
\item
This is followed by the eight combinations of Constraints B--D added to Constraint A. Along the line, the expected number of intra-continental matches can be reduced by almost two with increasing inequality by about 38\%. The Skip mechanism remains competitive with the Uniform.
\item
The eight sets of restrictions containing Constraint E without Constraint A are essentially dominated by the second curve as the reduction in the number of unattractive matches is minimal for a higher value of inequality. The Skip mechanism increases $I$ by approximately 0.011, which is roughly equivalent to the half of the original level of inequality, caused by the selection of the three hosts.
\item
The last curve is established by the eight combinations of Constraints B--D added to Constraints A and E. Even though this line provides the lowest number of intra-continental matches, the price is substantial, especially if the Skip mechanism is used. Fixing the number of unattractive matches at 12 means that inequality is more than doubled. Furthermore, imposing only Constraints A and E is inefficient as Constraints A--D ensure fewer unattractive games, however, the level of inequality is reduced by 14\% (Uniform mechanism) or even by 30\% (Skip mechanism).
\item
For both the Skip and the Uniform mechanisms, eight sets of draw constraints (1, 3, 5, 7, 9, 11, 13, and 17 in Table~\ref{Table3}) are dominated by the set of Constraints A--D (scenario 30) as the latter implies fewer unattractive games and a lower level of inequality.
\end{itemize}

\input{Table4_assignment_probability_heatmap}

Last but not least, Table~\ref{Table4} compares the assignment probabilities under the Skip and the Uniform mechanisms if---similar to basketball and football---all geographical constraints are imposed on the group draw. The two draw procedures are equivalent for matching the teams in the two strongest Pots 1 and 2. The Skip mechanism essentially provides equal treatment for Pots 1--3, where Constraint E still has no effect. The only exception is in the relation of Algeria and Austria/Croatia since one of these European teams should play in the group of Egypt due to the choice of the hosts, hence, they have a smaller chance to play against the only African team drawn from Pot 3. However, the Skip mechanism is strongly distorted for Switzerland, the only European team drawn from Pot 4. For instance, it plays against Egypt with a probability of 47.1\% and against one of the seven European teams drawn from Pot 1 with a probability of about 7.6\% in a fair draw but the corresponding values are 76.7\% and 3.3\%, respectively, if the Skip mechanism is used. Analogously, the opponent of Austria/Croatia from Pot 4 is chosen with a remarkably higher bias by the Skip mechanism.

Table~\ref{Table4} suggests that the Skip mechanism can provide equal treatment in a group draw for more pairs of pots even if it remains more distorted for the teams drawn at the end. According to our knowledge, this potential advantage has never been recognised before. The issue requires further investigation because it can provide an argument for preferring the Skip mechanism to a Uniform draw, which contradicts the standard definition of fairness in previous studies.

\section{Conclusions} \label{Sec5}

First in the literature, this paper has identified and explored an interesting trade-off in a group draw of a sports tournament between the number of games played by teams from the same geographic zone and the extent of violating equal treatment. A metric based on the Herfindahl--Hirschman index, a popular measure of market concentration and competitive balance, has been suggested to quantify the level of inequality. The trade-off is governed by the constraints imposed on the draw.

For the 2025 IHF Men's World Handball Championship draw, 64 different scenarios given by 32 reasonable sets of geographic restrictions and two draw procedures have been simulated 1 million times to compute the value of both objective functions. The Pareto frontier is shown to consist of roughly three lines for both the Skip and the Uniform mechanisms. Eight combinations of constraints turn out to be suboptimal. Similar calculations are worth conducting before the draw of major sports competitions to inform decision-makers about the effects of draw constraints.

Our research can be continued in several directions. We have examined only one particular tournament with a given distribution of teams in the pots. The recent FIBA World Cups and FIFA World Cups, where, in contrast to our case study in handball, the organisers have imposed all possible geographical restrictions, may also be investigated. The role of the order in which the pots are drawn has not been studied. Finally, the Skip mechanism is not necessarily worse than the Uniform mechanism with respect to equal treatment according to Table~\ref{Table4}.

\section*{Acknowledgements}

This paper could not have been written without the \emph{father} of the first author (also called \emph{L\'aszl\'o Csat\'o}), who has helped to code the simulations in Python. \\
%We are grateful to \emph{Dries Goossens}, \emph{Alex Krumer}, \emph{Frits C.~R.~Spieksma}, and \emph{Stephan Westphal} for useful advice. \\
%Eight anonymous colleagues and \emph{Ilia Tsetlin} have provided valuable comments and suggestions on earlier drafts. \\
The research was supported by the National Research, Development and Innovation Office under Grants FK 145838 and PD 146055, by the J\'anos Bolyai Research Scholarship of the Hungarian Academy of Sciences, and by the EK\"OP-24 University Research Scholarship Program of the Ministry for Culture and Innovation from the source of the National Research, Development and Innovation Fund. 

\bibliographystyle{apalike} 
\bibliography{All_references}

\end{document}

%% file: Figure1_unattractive_match_distribution.tex
\begin{figure}[t!]
\centering

\begin{tikzpicture}
%\selectcolormodel{gray}
\begin{axis}[
width = \textwidth, 
height = 0.6\textwidth,
xmajorgrids,
ymajorgrids,
%xbar stacked,
%bar width = 10pt,
xlabel = {Number of intra-continental games},
xlabel style = {align=center, font=\small},
%xtick = data,
%symbolic x coords = {Final,Semifinals,Quarterfinals,Round of 16,Round of 32,Group stage},
enlarge x limits = 0.1,
ylabel = {Probability in percentages (\%)},
ylabel style = {align=center, font=\small},
yticklabel style = {/pgf/number format/fixed,/pgf/number format/precision=5},
scaled y ticks = false,
ymin = 0,
ybar = 4pt,
bar width = 16pt,
]
% Imbalanced format
\addplot [blue, thick, pattern = crosshatch dots, pattern color = blue] coordinates{
(12,5.616)
(13,16.328)
(14,27.556)
(15,25.2497)
(16,15.0574)
(17,6.9365)
(18,2.4335)
(19,0.6903)
(20,0.1177)
(21,0.0149)
};
\end{axis}
\end{tikzpicture}

\captionsetup{justification=centerfirst}
\caption{The distribution of unattractive group matches without draw constraints}
\label{Fig1}

\end{figure}
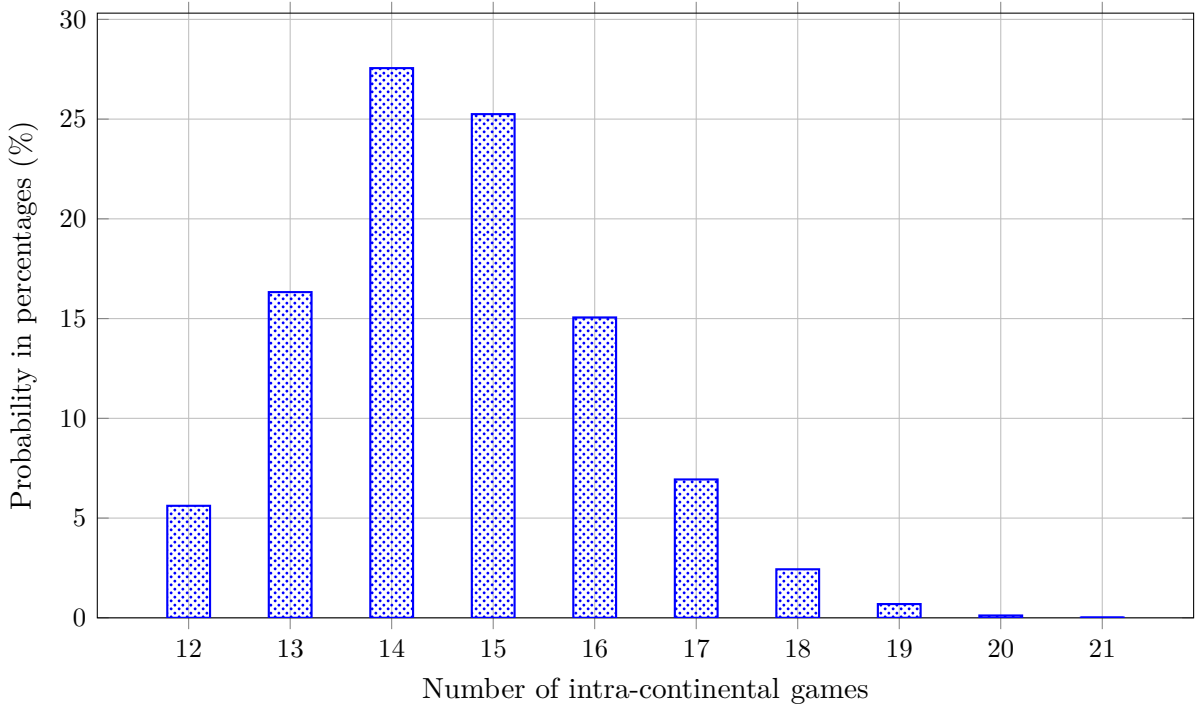

%\end{document}

%% file: Figure2_constraints_impact_feasibility.tex
\begin{figure}[t!]
\centering

\begin{tikzpicture}
\begin{axis}[
axis y line* = left,
xlabel = Set of constraints imposed (see Table~\ref{Table3}),
x label style = {font=\small},
ylabel = Proportion of valid assignments,
y label style = {font=\small},
%y tick label style = {/pgf/number format/.cd,fixed,fixed zerofill,precision=2},
width = 0.9\textwidth,
height = 0.6\textwidth,
ymajorgrids = true,
xmin = -0.5,
xmax = 31.5,
ymin = 0,
ymax = 1.5,
%max space between ticks=50,
%legend style = {font=\small,at={(0,-0.25)},anchor=north west,legend columns=5},
%legend entries = {15 states$\quad$,17 states$\quad$,20 states}
] 
% Ratio without Constraint E
\addplot [blue, mark=asterisk, only marks, mark size=2pt, mark options={solid,semithick}] coordinates {
(0,1)
(2,0.750461721574199)
(4,0.875491917020876)
(6,0.660358918279263)
(8,0.535881557314676)
(10,0.416506832169821)
(12,0.476512462468687)
(14,0.371527334009018)
(16,0.312366463336923)
(18,0.240750370755571)
(20,0.276448368291151)
(22,0.214451868529853)
(24,0.178597484668745)
(26,0.140456579393503)
(28,0.159737595856535)
(30,0.125972158641274)
};
% Ratio with Constraint E
\addplot [blue, mark=pentagon, only marks, mark size=2pt, mark options={solid,semithick}] coordinates {
(1,0.312795200470444)
(3,0.236595112418168)
(5,0.274096208865697)
(7,0.209016425137737)
(9,0.171185185223226)
(11,0.133660641512269)
(13,0.152360660843036)
(15,0.119215414743871)
(17,0.125020628403687)
(19,0.0998490282692569)
(21,0.112463813362253)
(23,0.0899764441669171)
(25,0.0774569339447267)
(27,0.0622278116625498)
(29,0.0698950735156383)
(31,0.0562087512079261)
};
\end{axis}

\begin{axis}[
axis y line* = right,
%title = {The effects of draw constraints},
%title style = {font=\small},
%xlabel = Set of constraints imposed (see Table),
%x label style = {font=\small},
xticklabel = \empty,
ylabel = {Proportion of the inequalities of the \\ Skip and Uniform mechanisms},
y label style = {font=\small,align=center},
%y tick label style = {/pgf/number format/.cd,fixed,fixed zerofill,precision=2},
width = 0.9\textwidth,
height = 0.6\textwidth,
ymajorgrids = true,
xmin = -0.5,
xmax = 31.5,
ymin = 0,
ymax = 1.5,
%max space between ticks=50,
%legend style = {font=\small,at={(0,-0.25)},anchor=north west,legend columns=5},
%legend entries = {11 states$\quad$,13 states$\quad$,15 states}
] 
% Relative inequality distortion of the Skip mechanism without Constraint E
\addplot [ForestGreen, mark=triangle, only marks, mark size=2pt, mark options={solid,semithick}] coordinates {
(0,0.999997946825151)
(2,0.999995580036349)
(4,1.00000246037314)
(6,0.999914116860847)
(8,1.00000734973822)
(10,0.999335732336326)
(12,0.999781234638578)
(14,0.999526788104473)
(16,1.0001646734137)
(18,1.00018533546766)
(20,0.999997280365302)
(22,1.00034567476295)
(24,1.00075475602635)
(26,1.00122735243373)
(28,1.00086018787915)
(30,1.00105552053626)
};
% Relative inequality distortion of the Skip mechanism with Constraint E but without Constraint A
\addplot [ForestGreen, mark=square, only marks, mark size=2pt, mark options={solid,thick}] coordinates {
(1,1.34151687387246)
(3,1.33121753929585)
(5,1.3389659067379)
(7,1.32690586243309)
(9,1.31690613844328)
(11,1.30959351116795)
(13,1.31367494606493)
(15,1.30813957824913)
};
% Relative inequality distortion of the Skip mechanism with Constraints A and E
\addplot [ForestGreen, mark=o, only marks, mark size=2pt, mark options={solid,semithick}] coordinates {
(17,1.2409730919088)
(19,1.23264563803348)
(21,1.23868348091692)
(23,1.23310332096554)
(25,1.22689028655375)
(27,1.2235537467907)
(29,1.22627587407787)
(31,1.2244435026637)
};
\end{axis}
\end{tikzpicture}

\begin{tikzpicture}
\begin{customlegend}[legend entries={Ratio of feasible solutions without Constraint E (left scale)$\qquad \qquad \qquad \quad \; \;$, Ratio of feasible solutions with Constraint E (left scale)$\qquad \qquad \qquad \qquad \quad$, Relative distortion without Constraint E (right scale)$\qquad \qquad \qquad \qquad \qquad \,$, Relative distortion with Constraint E and without Constraint A (right scale), Relative distortion with Constraints A and E (right scale)$\qquad \qquad \qquad \qquad \;$},legend columns=1]
    \addlegendimage{blue, mark=asterisk, only marks, mark size=2pt, mark options={solid,semithick}}
    \addlegendimage{blue, mark=pentagon, only marks, mark size=2pt, mark options={solid,semithick}}
    \addlegendimage{ForestGreen, mark=triangle, only marks, mark size=2pt, mark options={solid,semithick}}
    \addlegendimage{ForestGreen, mark=square, only marks, mark size=2pt, mark options={solid,thick}}
    \addlegendimage{ForestGreen, mark=o, only marks, mark size=2pt, mark options={solid,semithick}}
\end{customlegend}
\end{tikzpicture}

%\captionsetup{justification=centerfirst}
\caption{The effects of draw constraints and mechanisms: \\ the chance of feasibility and relative inequality distortion}
\label{Fig2}
\end{figure}
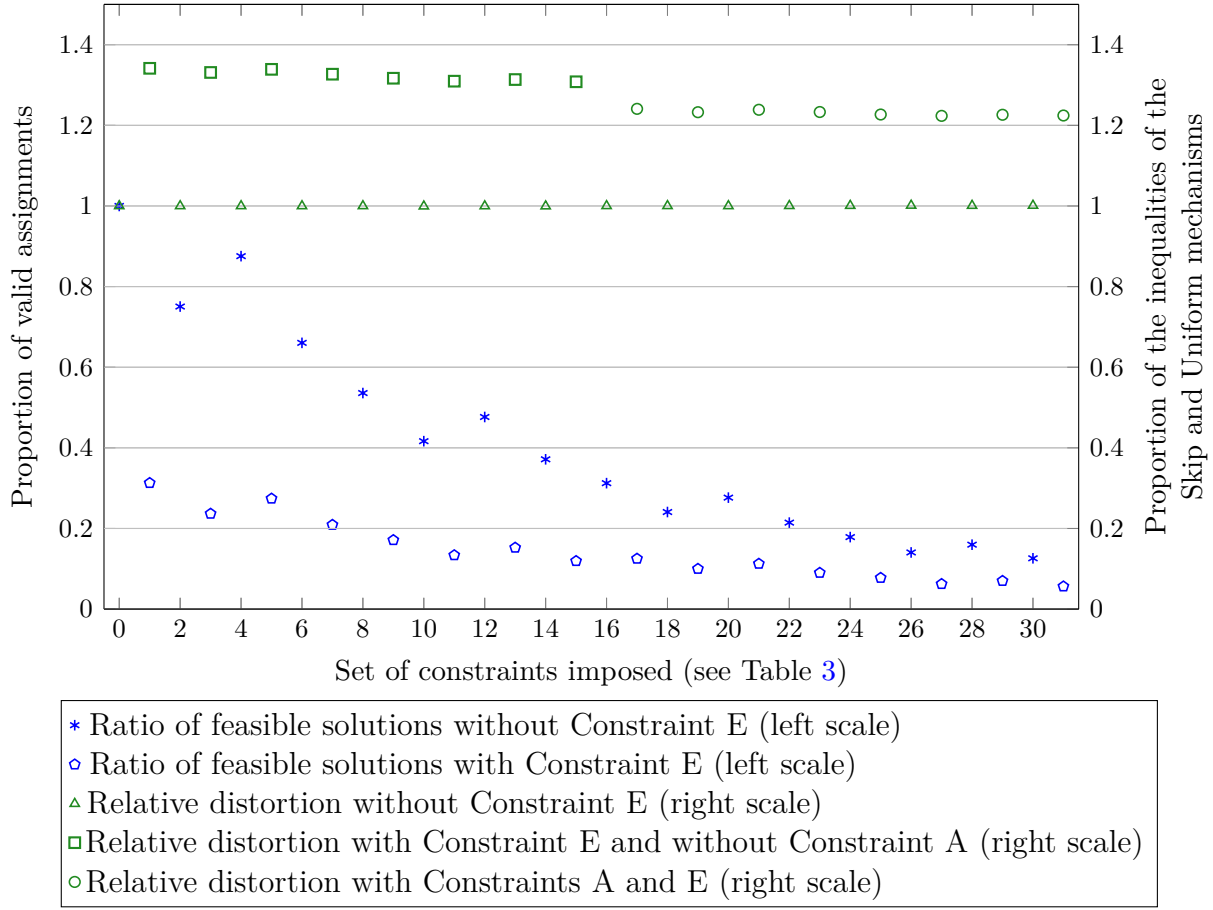

%\end{document}

%% file: Figure3_attractiveness_inequality_tradeoff.tex
\begin{figure}[t!]
\centering

\begin{tikzpicture}
\begin{axis}[
name = axis1,
xlabel = Average number of unattractive matches,
x label style = {font=\small},
ylabel = Level of inequality ($I$),
y label style = {font=\small},
yticklabel style = {/pgf/number format/fixed,/pgf/number format/precision=5},
scaled y ticks = false,
width = 0.96\textwidth,
height = 0.6\textwidth,
ymajorgrids = true,
xmin = 11.95,
xmax = 14.75,
ymin = 0.019,
ymax = 0.055,
%max space between ticks=50,
legend style = {font=\small,at={(0,-0.15)},anchor=north west,legend columns=1},
legend entries = {Uniform mechanism without Constraints A and E$\qquad \qquad \qquad \quad$, Uniform mechanism with Constraint A and without Constraint E, Uniform mechanism with Constraint E and without Constraint A, Uniform mechanism with Constraints A and E$\qquad \qquad \qquad \qquad \;$, Skip mechanism without Constraints A and E$\qquad \qquad \qquad \qquad \; \,$, Skip mechanism with Constraint A and without Constraint E$\quad \; \,$, Skip mechanism with Constraint E$\qquad \qquad \qquad \qquad \qquad \quad \qquad \;$}
] 
% Uniform without Constraints A and E
\addplot [blue, mark=asterisk, only marks, mark size=2pt, mark options={solid,semithick}] coordinates {
(14.623735,0.0238104880571905)
(14.314579,0.0249442258220476)
(14.472506,0.0242963920630952)
(14.17893,0.0253876917887619)
(14.013738,0.0264559276742381)
(13.777432,0.0273755037301905)
(13.895386,0.0268467486282857)
(13.660895,0.0277231062673333)
};
% Uniform with Constraint A and without Constraint E
\addplot [blue, mark=10-pointed star, only marks, mark size=2pt, mark options={solid,semithick}] coordinates {
(13.54183,0.0294584048892857)
(13.278539,0.0304680620229047)
(13.412804,0.0298854562648571)
(13.154347,0.0308294293046191)
(13.014954,0.0318583463600952)
(12.785918,0.0325705640831429)
(12.900646,0.0321537861268571)
(12.677143,0.0328384972601905)
};
% Uniform with Constraint E and without Constraint A
\addplot [blue, mark=pentagon, only marks, mark size=2pt, mark options={solid,semithick}] coordinates {
(13.521361,0.0328393191032381)
(13.223371,0.0338423472790476)
(13.374402,0.0332345826533333)
(13.090548,0.0342088881940476)
(12.929201,0.0351894029060952)
(12.694118,0.0359823891620952)
(12.812659,0.0355283881432857)
(12.580365,0.0362714296960952)
};
% Uniform with Constraints A and E
\addplot [blue, mark=Mercedes star, only marks, mark size=2pt, mark options={solid,semithick}] coordinates {
(12.701015,0.0379450454715714)
(12.498135,0.0387446924452381)
(12.598717,0.0382795040561905)
(12.396396,0.0389892283243333)
(12.296542,0.0397658431159048)
(12.096938,0.0402742189039524)
(12.195412,0.0399601546144286)
(12,0.0404590970337143)
};
% Skip mechanism with Constraint A and without Constraint E
\addplot [red, mark=o, only marks, mark size=2pt, mark options={solid,thin}] coordinates {
(14.625858,0.0238104391700952)
(14.315941,0.0249441155694762)
(14.47444,0.0242964518412857)
(14.175213,0.0253855114140952)
(14.014972,0.026456122118381)
(13.767612,0.0273573190682857)
(13.888548,0.0268408754896191)
(13.648896,0.0277099873636667)
};
% Skip mechanism without Constraint E
\addplot [red, mark=square, only marks, mark size=2pt, mark options={solid,thin}] coordinates {
(13.544032,0.0294632559053809)
(13.271508,0.0304737088354286)
(13.407855,0.0298853749873333)
(13.144999,0.0308400862602857)
(12.999586,0.031882391639)
(12.778397,0.0326105396442381)
(12.887823,0.0321814444239524)
(12.673859,0.0328731589684286)
};
% Skip mechanism with Constraint E
\addplot [red, mark=oplus, only marks, mark size=2pt, mark options={solid,thin}] coordinates {
(13.416706,0.0440545007034762)
(13.123415,0.0450515262688095)
(13.271256,0.0444999730974762)
(12.993344,0.045391974292)
(12.84073,0.0463411406951905)
(12.609909,0.047122303363)
(12.725246,0.0466727533779048)
(12.502446,0.0474480927451429)
(12.697531,0.0470887804014762)
(12.480997,0.0477584761395714)
(12.588286,0.0474161893320952)
(12.37884,0.0480777469286191)
(12.273829,0.0487883266555238)
(12.087942,0.049277671439)
(12.180281,0.0490021735280952)
(12,0.0495398784865714)
};
\end{axis}
\end{tikzpicture}

%\captionsetup{justification=centerfirst}
\caption{The trade-off between unattractive matches and the inequality of the draw}
\label{Fig3} 
\end{figure}
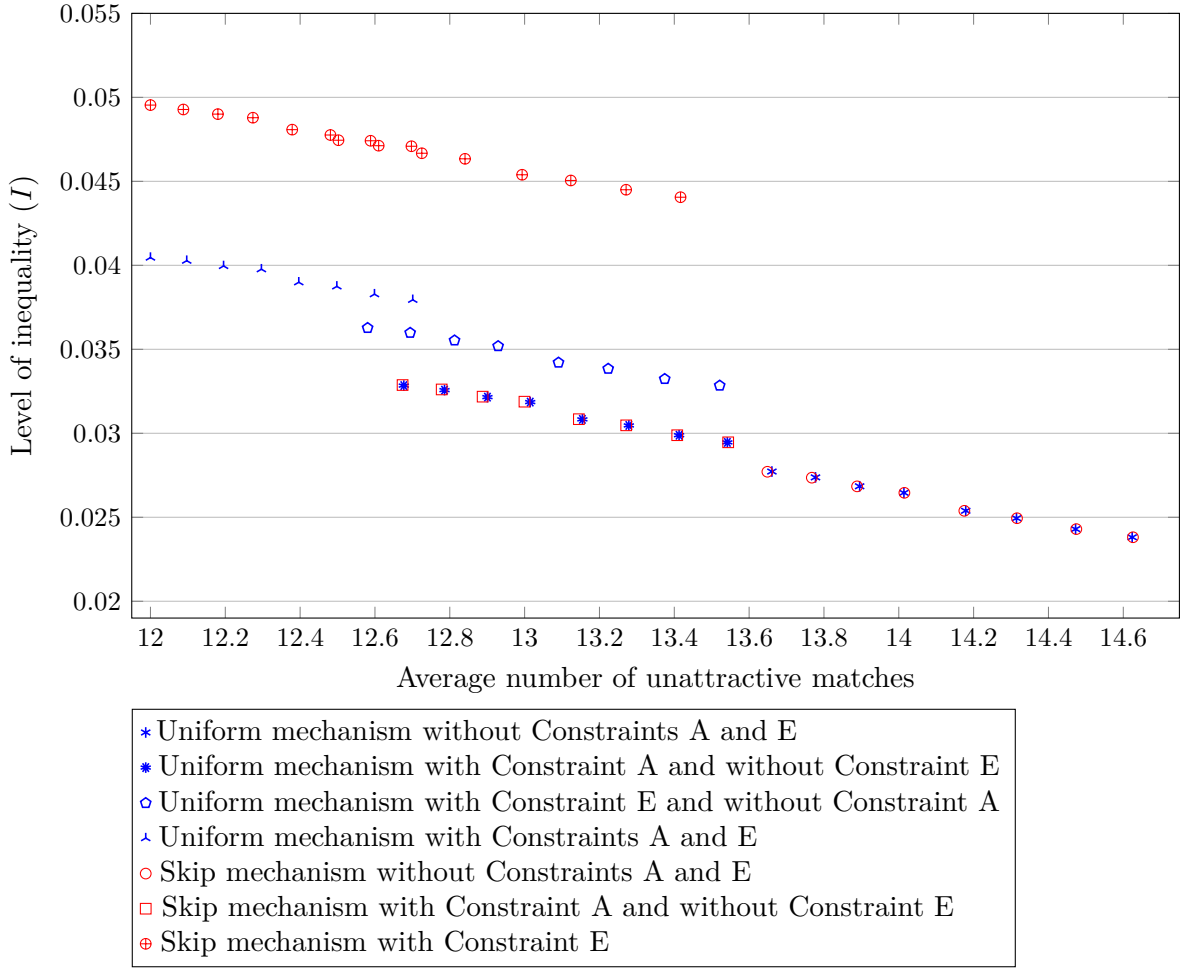

%\end{document}

%% file: Table4_assignment_probability_heatmap.tex
\begin{table}[t!]
  \centering
  \caption{Assignment probabilities in the 2025 IHF World Men's Handball \\ Championship draw with all geographical Constraints A--E}
  \label{Table4}
  
\begin{subtable}{\textwidth}
  \caption{Uniform mechanism (draw order 1-2-3-4)}
  \label{Table4a}
\resizebox{1\textwidth}{!}{
\begin{tiny}
    \begin{tabularx}{1.4\textwidth}{r CCCC CCCC CCCC CCCC CCCC CCCC} \toprule
          & \rotatebox[origin=l]{90}{Portugal} & \rotatebox[origin=l]{90}{Croatia} & \rotatebox[origin=l]{90}{Austria} & \rotatebox[origin=l]{90}{Iceland} & \rotatebox[origin=l]{90}{Netherlands} & \rotatebox[origin=l]{90}{Spain} & \rotatebox[origin=l]{90}{Italy} & \rotatebox[origin=l]{90}{Czechia} & \rotatebox[origin=l]{90}{Poland} & \rotatebox[origin=l]{90}{Macedonia} & \rotatebox[origin=l]{90}{Qatar} & \rotatebox[origin=l]{90}{Brazil} & \rotatebox[origin=l]{90}{Argentina} & \rotatebox[origin=l]{90}{Cuba} & \rotatebox[origin=l]{90}{Japan} & \rotatebox[origin=l]{90}{Algeria} & \rotatebox[origin=l]{90}{Bahrain} & \rotatebox[origin=l]{90}{Tunisia} & \rotatebox[origin=l]{90}{Chile} & \rotatebox[origin=l]{90}{Kuwait} & \rotatebox[origin=l]{90}{Cape Verde} & \rotatebox[origin=l]{90}{Guinea} & \rotatebox[origin=l]{90}{United States} & \rotatebox[origin=l]{90}{Switzerland} \\ \bottomrule
    Denmark & \cellcolor{ForestGreen!4.77931428571429} 16.7 & X     & X     & \cellcolor{ForestGreen!4.70468571428571} 16.6 & \cellcolor{ForestGreen!4.78857142857143} 16.7 & \cellcolor{ForestGreen!4.7824} 16.7 & \cellcolor{ForestGreen!4.75417142857143} 16.7 & \cellcolor{ForestGreen!4.76228571428571} 16.7 & \cellcolor{red!22.0296} 9.7 & \cellcolor{red!22.384} 9.7 & \cellcolor{ForestGreen!0.733028571428571} 13.1 & \cellcolor{ForestGreen!0.889371428571429} 13.3 & \cellcolor{ForestGreen!0.931771428571429} 13.3 & \cellcolor{ForestGreen!1.04045714285714} 13.4 & \cellcolor{ForestGreen!0.751771428571429} 13.2 & \cellcolor{ForestGreen!1.9984} 14.2 & \cellcolor{red!2.1968} 12.2 & \cellcolor{ForestGreen!2.04091428571429} 14.3 & \cellcolor{red!0.6072} 12.4 & \cellcolor{red!1.7544} 12.3 & \cellcolor{ForestGreen!2.02868571428571} 14.3 & \cellcolor{ForestGreen!2.07954285714286} 14.3 & \cellcolor{ForestGreen!0.195885714285714} 12.7 & \cellcolor{red!39.8568} 7.5 \\
    France & X     & \cellcolor{ForestGreen!42.8990857142857} 50 & \cellcolor{ForestGreen!42.8152} 50 & X     & X     & X     & X     & X     & \cellcolor{red!22.8368} 9.6 & \cellcolor{red!22.4152} 9.7 & \cellcolor{ForestGreen!0.683657142857143} 13.1 & \cellcolor{ForestGreen!0.991428571428571} 13.4 & \cellcolor{ForestGreen!1.00571428571429} 13.4 & \cellcolor{ForestGreen!0.989714285714286} 13.4 & \cellcolor{ForestGreen!0.742171428571429} 13.1 & \cellcolor{ForestGreen!2.05188571428571} 14.3 & \cellcolor{red!1.6432} 12.3 & \cellcolor{ForestGreen!2.06034285714286} 14.3 & \cellcolor{red!0.9048} 12.4 & \cellcolor{red!1.6416} 12.3 & \cellcolor{ForestGreen!2.06685714285714} 14.3 & \cellcolor{ForestGreen!1.96742857142857} 14.2 & \cellcolor{ForestGreen!0.162171428571429} 12.6 & \cellcolor{red!39.608} 7.5 \\
    Sweden & \cellcolor{ForestGreen!4.68902857142857} 16.6 & X     & X     & \cellcolor{ForestGreen!4.84262857142857} 16.7 & \cellcolor{ForestGreen!4.73794285714286} 16.6 & \cellcolor{ForestGreen!4.76514285714286} 16.7 & \cellcolor{ForestGreen!4.76445714285714} 16.7 & \cellcolor{ForestGreen!4.77222857142857} 16.7 & \cellcolor{red!22.632} 9.7 & \cellcolor{red!22.3488} 9.7 & \cellcolor{ForestGreen!0.735085714285714} 13.1 & \cellcolor{ForestGreen!1.03074285714286} 13.4 & \cellcolor{ForestGreen!0.988114285714286} 13.4 & \cellcolor{ForestGreen!0.950514285714286} 13.3 & \cellcolor{ForestGreen!0.748457142857143} 13.2 & \cellcolor{ForestGreen!1.97291428571429} 14.2 & \cellcolor{red!1.4992} 12.3 & \cellcolor{ForestGreen!2.09611428571429} 14.3 & \cellcolor{red!1.0152} 12.4 & \cellcolor{red!1.5536} 12.3 & \cellcolor{ForestGreen!1.98834285714286} 14.2 & \cellcolor{ForestGreen!2.05382857142857} 14.3 & \cellcolor{ForestGreen!0.104571428571429} 12.6 & \cellcolor{red!39.632} 7.5 \\
    Germany & \cellcolor{ForestGreen!4.79634285714286} 16.7 & X     & X     & \cellcolor{ForestGreen!4.75302857142857} 16.7 & \cellcolor{ForestGreen!4.77451428571429} 16.7 & \cellcolor{ForestGreen!4.78628571428571} 16.7 & \cellcolor{ForestGreen!4.74502857142857} 16.7 & \cellcolor{ForestGreen!4.71622857142857} 16.6 & \cellcolor{red!22.2392} 9.7 & \cellcolor{red!22.528} 9.7 & \cellcolor{ForestGreen!0.737485714285714} 13.1 & \cellcolor{ForestGreen!0.932} 13.3 & \cellcolor{ForestGreen!0.9336} 13.3 & \cellcolor{ForestGreen!0.954742857142857} 13.3 & \cellcolor{ForestGreen!0.779885714285714} 13.2 & \cellcolor{ForestGreen!2.0576} 14.3 & \cellcolor{red!1.9008} 12.3 & \cellcolor{ForestGreen!2.02285714285714} 14.3 & \cellcolor{red!0.9656} 12.4 & \cellcolor{red!1.9888} 12.3 & \cellcolor{ForestGreen!2.11302857142857} 14.3 & \cellcolor{ForestGreen!2.01622857142857} 14.3 & \cellcolor{ForestGreen!0.171885714285714} 12.7 & \cellcolor{red!39.4128} 7.6 \\
    Hungary & \cellcolor{ForestGreen!4.77188571428571} 16.7 & X     & X     & \cellcolor{ForestGreen!4.79851428571429} 16.7 & \cellcolor{ForestGreen!4.79314285714286} 16.7 & \cellcolor{ForestGreen!4.76708571428571} 16.7 & \cellcolor{ForestGreen!4.69474285714286} 16.6 & \cellcolor{ForestGreen!4.74605714285714} 16.7 & \cellcolor{red!22.6464} 9.7 & \cellcolor{red!22.4} 9.7 & \cellcolor{ForestGreen!0.793142857142857} 13.2 & \cellcolor{ForestGreen!0.937028571428571} 13.3 & \cellcolor{ForestGreen!0.909142857142857} 13.3 & \cellcolor{ForestGreen!1.012} 13.4 & \cellcolor{ForestGreen!0.726057142857143} 13.1 & \cellcolor{ForestGreen!2.05782857142857} 14.3 & \cellcolor{red!1.8776} 12.3 & \cellcolor{ForestGreen!2.02434285714286} 14.3 & \cellcolor{red!0.9992} 12.4 & \cellcolor{red!1.728} 12.3 & \cellcolor{ForestGreen!2.06045714285714} 14.3 & \cellcolor{ForestGreen!2.0576} 14.3 & \cellcolor{ForestGreen!0.175542857142857} 12.7 & \cellcolor{red!39.6208} 7.5 \\
    Slovenia & \cellcolor{ForestGreen!4.77817142857143} 16.7 & X     & X     & \cellcolor{ForestGreen!4.7464} 16.7 & \cellcolor{ForestGreen!4.75554285714286} 16.7 & \cellcolor{ForestGreen!4.72057142857143} 16.6 & \cellcolor{ForestGreen!4.78045714285714} 16.7 & \cellcolor{ForestGreen!4.79028571428571} 16.7 & \cellcolor{red!22.872} 9.6 & \cellcolor{red!22.2056} 9.7 & \cellcolor{ForestGreen!0.799428571428571} 13.2 & \cellcolor{ForestGreen!0.9288} 13.3 & \cellcolor{ForestGreen!0.9744} 13.4 & \cellcolor{ForestGreen!0.971885714285714} 13.4 & \cellcolor{ForestGreen!0.716914285714286} 13.1 & \cellcolor{ForestGreen!2.04822857142857} 14.3 & \cellcolor{red!1.5144} 12.3 & \cellcolor{ForestGreen!2.04708571428571} 14.3 & \cellcolor{red!0.8728} 12.4 & \cellcolor{red!1.9696} 12.3 & \cellcolor{ForestGreen!2.03874285714286} 14.3 & \cellcolor{ForestGreen!2.04628571428571} 14.3 & \cellcolor{ForestGreen!0.121371428571429} 12.6 & \cellcolor{red!39.4176} 7.6 \\
    Norway & \cellcolor{ForestGreen!4.75668571428571} 16.7 & X     & X     & \cellcolor{ForestGreen!4.72617142857143} 16.6 & \cellcolor{ForestGreen!4.72171428571429} 16.6 & \cellcolor{ForestGreen!4.74994285714286} 16.7 & \cellcolor{ForestGreen!4.83257142857143} 16.7 & \cellcolor{ForestGreen!4.78434285714286} 16.7 & \cellcolor{red!22.9296} 9.6 & \cellcolor{red!22.732} 9.7 & \cellcolor{ForestGreen!0.756685714285714} 13.2 & \cellcolor{ForestGreen!0.9896} 13.4 & \cellcolor{ForestGreen!0.987428571428571} 13.4 & \cellcolor{ForestGreen!0.937371428571429} 13.3 & \cellcolor{ForestGreen!0.753142857142857} 13.2 & \cellcolor{ForestGreen!2.09885714285714} 14.3 & \cellcolor{red!1.792} 12.3 & \cellcolor{ForestGreen!1.99405714285714} 14.2 & \cellcolor{red!0.6072} 12.4 & \cellcolor{red!1.6656} 12.3 & \cellcolor{ForestGreen!1.9896} 14.2 & \cellcolor{ForestGreen!2.0648} 14.3 & \cellcolor{ForestGreen!0.169028571428571} 12.6 & \cellcolor{red!39.4576} 7.6 \\
    Egypt & X     & \cellcolor{ForestGreen!42.8152} 50 & \cellcolor{ForestGreen!42.8990857142857} 50 & X     & X     & X     & X     & X     & \cellcolor{ForestGreen!22.5979428571429} 32.3 & \cellcolor{ForestGreen!22.4305142857143} 32.1 & \cellcolor{red!36.6696} 7.9 & \cellcolor{red!46.8928} 6.6 & \cellcolor{red!47.1112} 6.6 & \cellcolor{red!47.9968} 6.5 & \cellcolor{red!36.5288} 7.9 & X     & \cellcolor{ForestGreen!1.77485714285714} 14.1 & X     & \cellcolor{ForestGreen!0.853142857142857} 13.2 & \cellcolor{ForestGreen!1.75737142857143} 14 & X     & X     & \cellcolor{red!7.7032} 11.5 & \cellcolor{ForestGreen!39.5722285714286} 47.1 \\
    Portugal & X     & X     & X     & X     & X     & X     & X     & X     & \cellcolor{red!22.9672} 9.6 & \cellcolor{red!22.5136} 9.7 & \cellcolor{ForestGreen!0.748114285714286} 13.2 & \cellcolor{ForestGreen!0.927314285714286} 13.3 & \cellcolor{ForestGreen!0.939885714285714} 13.3 & \cellcolor{ForestGreen!1.03702857142857} 13.4 & \cellcolor{ForestGreen!0.767085714285714} 13.2 & \cellcolor{ForestGreen!2.07782857142857} 14.3 & \cellcolor{red!2.208} 12.2 & \cellcolor{ForestGreen!1.97085714285714} 14.2 & \cellcolor{red!0.8} 12.4 & \cellcolor{red!1.0656} 12.4 & \cellcolor{ForestGreen!2.09508571428571} 14.3 & \cellcolor{ForestGreen!2.03577142857143} 14.3 & \cellcolor{ForestGreen!0.165714285714286} 12.6 & \cellcolor{red!39.7984} 7.5 \\
    Croatia & X     & X     & X     & X     & X     & X     & X     & X     & \cellcolor{ForestGreen!9.68708571428571} 21 & \cellcolor{ForestGreen!9.56982857142857} 20.9 & \cellcolor{red!16.0096} 10.5 & \cellcolor{red!20.0176} 10 & \cellcolor{red!19.6648} 10 & \cellcolor{red!20.2808} 10 & \cellcolor{red!15.7776} 10.5 & \cellcolor{red!43.048} 7.1 & \cellcolor{ForestGreen!0.733485714285714} 13.1 & \cellcolor{red!42.9808} 7.1 & \cellcolor{ForestGreen!0.365828571428571} 12.8 & \cellcolor{ForestGreen!0.792457142857143} 13.2 & \cellcolor{red!42.4976} 7.2 & \cellcolor{red!42.82} 7.1 & \cellcolor{red!3.392} 12.1 & \cellcolor{ForestGreen!16.9211428571429} 27.3 \\
    Austria & X     & X     & X     & X     & X     & X     & X     & X     & \cellcolor{ForestGreen!9.64845714285714} 20.9 & \cellcolor{ForestGreen!9.65851428571429} 21 & \cellcolor{red!15.8744} 10.5 & \cellcolor{red!19.9352} 10 & \cellcolor{red!20.4064} 9.9 & \cellcolor{red!20.788} 9.9 & \cellcolor{red!15.556} 10.6 & \cellcolor{red!42.5888} 7.2 & \cellcolor{ForestGreen!0.806628571428571} 13.2 & \cellcolor{red!42.5968} 7.2 & \cellcolor{ForestGreen!0.358057142857143} 12.8 & \cellcolor{ForestGreen!0.7304} 13.1 & \cellcolor{red!43.0344} 7.1 & \cellcolor{red!43.408} 7.1 & \cellcolor{red!3.176} 12.1 & \cellcolor{ForestGreen!16.9928} 27.4 \\
    Iceland & X     & X     & X     & X     & X     & X     & X     & X     & \cellcolor{red!22.6728} 9.7 & \cellcolor{red!22.9368} 9.6 & \cellcolor{ForestGreen!0.78} 13.2 & \cellcolor{ForestGreen!0.972457142857143} 13.4 & \cellcolor{ForestGreen!1.01085714285714} 13.4 & \cellcolor{ForestGreen!1.04582857142857} 13.4 & \cellcolor{ForestGreen!0.676} 13.1 & \cellcolor{ForestGreen!2.03051428571429} 14.3 & \cellcolor{red!1.5296} 12.3 & \cellcolor{ForestGreen!2.0584} 14.3 & \cellcolor{red!0.9088} 12.4 & \cellcolor{red!2.2456} 12.2 & \cellcolor{ForestGreen!2.02022857142857} 14.3 & \cellcolor{ForestGreen!2.07577142857143} 14.3 & \cellcolor{ForestGreen!0.136342857142857} 12.6 & \cellcolor{red!39.3512} 7.6 \\
    Netherlands & X     & X     & X     & X     & X     & X     & X     & X     & \cellcolor{red!22.6104} 9.7 & \cellcolor{red!22.1544} 9.7 & \cellcolor{ForestGreen!0.748914285714286} 13.2 & \cellcolor{ForestGreen!0.9496} 13.3 & \cellcolor{ForestGreen!0.9736} 13.4 & \cellcolor{ForestGreen!0.937142857142857} 13.3 & \cellcolor{ForestGreen!0.767542857142857} 13.2 & \cellcolor{ForestGreen!2.01817142857143} 14.3 & \cellcolor{red!1.7216} 12.3 & \cellcolor{ForestGreen!2.06537142857143} 14.3 & \cellcolor{red!0.8544} 12.4 & \cellcolor{red!2.0088} 12.2 & \cellcolor{ForestGreen!2.07542857142857} 14.3 & \cellcolor{ForestGreen!1.99497142857143} 14.2 & \cellcolor{ForestGreen!0.190742857142857} 12.7 & \cellcolor{red!39.7008} 7.5 \\
    Spain & X     & X     & X     & X     & X     & X     & X     & X     & \cellcolor{red!21.9048} 9.8 & \cellcolor{red!22.8856} 9.6 & \cellcolor{ForestGreen!0.777828571428571} 13.2 & \cellcolor{ForestGreen!0.924114285714286} 13.3 & \cellcolor{ForestGreen!0.930971428571429} 13.3 & \cellcolor{ForestGreen!0.962971428571429} 13.3 & \cellcolor{ForestGreen!0.743085714285714} 13.2 & \cellcolor{ForestGreen!2.05965714285714} 14.3 & \cellcolor{red!1.6864} 12.3 & \cellcolor{ForestGreen!2.05874285714286} 14.3 & \cellcolor{red!0.944} 12.4 & \cellcolor{red!1.4672} 12.3 & \cellcolor{ForestGreen!1.96182857142857} 14.2 & \cellcolor{ForestGreen!2.03897142857143} 14.3 & \cellcolor{ForestGreen!0.156685714285714} 12.6 & \cellcolor{red!39.416} 7.6 \\
    Italy & X     & X     & X     & X     & X     & X     & X     & X     & \cellcolor{red!22.624} 9.7 & \cellcolor{red!21.9288} 9.8 & \cellcolor{ForestGreen!0.768342857142857} 13.2 & \cellcolor{ForestGreen!0.9752} 13.4 & \cellcolor{ForestGreen!0.877485714285714} 13.3 & \cellcolor{ForestGreen!0.905257142857143} 13.3 & \cellcolor{ForestGreen!0.789028571428571} 13.2 & \cellcolor{ForestGreen!2.04937142857143} 14.3 & \cellcolor{red!1.9568} 12.3 & \cellcolor{ForestGreen!2.07234285714286} 14.3 & \cellcolor{red!1.1488} 12.4 & \cellcolor{red!1.8128} 12.3 & \cellcolor{ForestGreen!2.00217142857143} 14.3 & \cellcolor{ForestGreen!2.1096} 14.3 & \cellcolor{ForestGreen!0.154971428571429} 12.6 & \cellcolor{red!39.4552} 7.6 \\
    Czechia & X     & X     & X     & X     & X     & X     & X     & X     & \cellcolor{red!22.5696} 9.7 & \cellcolor{red!22.1792} 9.7 & \cellcolor{ForestGreen!0.731657142857143} 13.1 & \cellcolor{ForestGreen!0.958857142857143} 13.3 & \cellcolor{ForestGreen!0.991657142857143} 13.4 & \cellcolor{ForestGreen!0.978742857142857} 13.4 & \cellcolor{ForestGreen!0.733485714285714} 13.1 & \cellcolor{ForestGreen!1.99828571428571} 14.2 & \cellcolor{red!1.6784} 12.3 & \cellcolor{ForestGreen!1.99965714285714} 14.2 & \cellcolor{red!0.4112} 12.4 & \cellcolor{red!2.06} 12.2 & \cellcolor{ForestGreen!2.06411428571429} 14.3 & \cellcolor{ForestGreen!2.0632} 14.3 & \cellcolor{ForestGreen!0.133828571428571} 12.6 & \cellcolor{red!39.676} 7.5 \\
    Poland & \cellcolor{red!22.9672} 9.6 & \cellcolor{ForestGreen!9.68708571428571} 21 & \cellcolor{ForestGreen!9.64845714285714} 20.9 & \cellcolor{red!22.6728} 9.7 & \cellcolor{red!22.6104} 9.7 & \cellcolor{red!21.9048} 9.8 & \cellcolor{red!22.624} 9.7 & \cellcolor{red!22.5696} 9.7 & X     & X     & X     & X     & X     & X     & X     & X     & \cellcolor{ForestGreen!5.09634285714286} 17 & \cellcolor{red!20.5344} 9.9 & \cellcolor{ForestGreen!4.2984} 16.3 & \cellcolor{ForestGreen!5.11531428571429} 17 & \cellcolor{red!20.188} 10 & \cellcolor{red!20.3008} 10 & \cellcolor{ForestGreen!1.9112} 14.2 & \cellcolor{red!53.9256} 5.8 \\
    Macedonia & \cellcolor{red!22.5136} 9.7 & \cellcolor{ForestGreen!9.56982857142857} 20.9 & \cellcolor{ForestGreen!9.65851428571429} 21 & \cellcolor{red!22.9368} 9.6 & \cellcolor{red!22.1544} 9.7 & \cellcolor{red!22.8856} 9.6 & \cellcolor{red!21.9288} 9.8 & \cellcolor{red!22.1792} 9.7 & X     & X     & X     & X     & X     & X     & X     & X     & \cellcolor{ForestGreen!5.06308571428571} 16.9 & \cellcolor{red!20.2256} 10 & \cellcolor{ForestGreen!4.22422857142857} 16.2 & \cellcolor{ForestGreen!5.09634285714286} 17 & \cellcolor{red!20.0152} 10 & \cellcolor{red!20.0112} 10 & \cellcolor{ForestGreen!1.91942857142857} 14.2 & \cellcolor{red!53.8696} 5.8 \\
    Qatar & \cellcolor{ForestGreen!0.748114285714286} 13.2 & \cellcolor{red!16.0096} 10.5 & \cellcolor{red!15.8744} 10.5 & \cellcolor{ForestGreen!0.78} 13.2 & \cellcolor{ForestGreen!0.748914285714286} 13.2 & \cellcolor{ForestGreen!0.777828571428571} 13.2 & \cellcolor{ForestGreen!0.768342857142857} 13.2 & \cellcolor{ForestGreen!0.731657142857143} 13.1 & X     & X     & X     & X     & X     & X     & X     & X     & X     & \cellcolor{ForestGreen!6.01577142857143} 17.8 & \cellcolor{ForestGreen!3.88422857142857} 15.9 & X     & \cellcolor{ForestGreen!5.97394285714286} 17.7 & \cellcolor{ForestGreen!5.93314285714286} 17.7 & \cellcolor{ForestGreen!1.75771428571429} 14 & \cellcolor{ForestGreen!5.00662857142857} 16.9 \\
    Brazil & \cellcolor{ForestGreen!0.927314285714286} 13.3 & \cellcolor{red!20.0176} 10 & \cellcolor{red!19.9352} 10 & \cellcolor{ForestGreen!0.972457142857143} 13.4 & \cellcolor{ForestGreen!0.9496} 13.3 & \cellcolor{ForestGreen!0.924114285714286} 13.3 & \cellcolor{ForestGreen!0.9752} 13.4 & \cellcolor{ForestGreen!0.958857142857143} 13.3 & X     & X     & X     & X     & X     & X     & X     & X     & \cellcolor{ForestGreen!2.01485714285714} 14.3 & \cellcolor{ForestGreen!2.81188571428571} 15 & X     & \cellcolor{ForestGreen!1.9736} 14.2 & \cellcolor{ForestGreen!2.77817142857143} 14.9 & \cellcolor{ForestGreen!2.82137142857143} 15 & \cellcolor{red!2.3512} 12.2 & \cellcolor{ForestGreen!2.22171428571429} 14.4 \\
    Argentina & \cellcolor{ForestGreen!0.939885714285714} 13.3 & \cellcolor{red!19.6648} 10 & \cellcolor{red!20.4064} 9.9 & \cellcolor{ForestGreen!1.01085714285714} 13.4 & \cellcolor{ForestGreen!0.9736} 13.4 & \cellcolor{ForestGreen!0.930971428571429} 13.3 & \cellcolor{ForestGreen!0.877485714285714} 13.3 & \cellcolor{ForestGreen!0.991657142857143} 13.4 & X     & X     & X     & X     & X     & X     & X     & X     & \cellcolor{ForestGreen!1.9584} 14.2 & \cellcolor{ForestGreen!2.80548571428571} 15 & X     & \cellcolor{ForestGreen!1.94468571428571} 14.2 & \cellcolor{ForestGreen!2.8392} 15 & \cellcolor{ForestGreen!2.84205714285714} 15 & \cellcolor{red!2.4568} 12.2 & \cellcolor{ForestGreen!2.24685714285714} 14.5 \\
    Cuba  & \cellcolor{ForestGreen!1.03702857142857} 13.4 & \cellcolor{red!20.2808} 10 & \cellcolor{red!20.788} 9.9 & \cellcolor{ForestGreen!1.04582857142857} 13.4 & \cellcolor{ForestGreen!0.937142857142857} 13.3 & \cellcolor{ForestGreen!0.962971428571429} 13.3 & \cellcolor{ForestGreen!0.905257142857143} 13.3 & \cellcolor{ForestGreen!0.978742857142857} 13.4 & X     & X     & X     & X     & X     & X     & X     & X     & \cellcolor{ForestGreen!1.79234285714286} 14.1 & \cellcolor{ForestGreen!2.47234285714286} 14.7 & \cellcolor{ForestGreen!1.31314285714286} 13.6 & \cellcolor{ForestGreen!1.66651428571429} 14 & \cellcolor{ForestGreen!2.51965714285714} 14.7 & \cellcolor{ForestGreen!2.5816} 14.8 & X     & \cellcolor{ForestGreen!1.94011428571429} 14.2 \\
    Japan & \cellcolor{ForestGreen!0.767085714285714} 13.2 & \cellcolor{red!15.7776} 10.5 & \cellcolor{red!15.556} 10.6 & \cellcolor{ForestGreen!0.676} 13.1 & \cellcolor{ForestGreen!0.767542857142857} 13.2 & \cellcolor{ForestGreen!0.743085714285714} 13.2 & \cellcolor{ForestGreen!0.789028571428571} 13.2 & \cellcolor{ForestGreen!0.733485714285714} 13.1 & X     & X     & X     & X     & X     & X     & X     & X     & X     & \cellcolor{ForestGreen!6.00308571428571} 17.8 & \cellcolor{ForestGreen!3.888} 15.9 & X     & \cellcolor{ForestGreen!5.91805714285714} 17.7 & \cellcolor{ForestGreen!5.8664} 17.6 & \cellcolor{ForestGreen!1.80422857142857} 14.1 & \cellcolor{ForestGreen!5.09165714285714} 17 \\
    Algeria & \cellcolor{ForestGreen!2.07782857142857} 14.3 & \cellcolor{red!43.048} 7.1 & \cellcolor{red!42.5888} 7.2 & \cellcolor{ForestGreen!2.03051428571429} 14.3 & \cellcolor{ForestGreen!2.01817142857143} 14.3 & \cellcolor{ForestGreen!2.05965714285714} 14.3 & \cellcolor{ForestGreen!2.04937142857143} 14.3 & \cellcolor{ForestGreen!1.99828571428571} 14.2 & X     & X     & X     & X     & X     & X     & X     & X     & \cellcolor{ForestGreen!12.6464} 23.6 & X     & \cellcolor{ForestGreen!10.9634285714286} 22.1 & \cellcolor{ForestGreen!12.7749714285714} 23.7 & X     & X     & \cellcolor{ForestGreen!7.58} 19.1 & \cellcolor{red!7.7536} 11.5 \\  \toprule
    \end{tabularx}
\end{tiny}
}
\end{subtable}

\vspace{0.25cm}
\begin{subtable}{\textwidth}
  \caption{Skip mechanism (draw order 1-2-3-4)}
  \label{Table4b}
\resizebox{\textwidth}{!}{
\begin{threeparttable}
\begin{tiny}
    \begin{tabularx}{1.4\textwidth}{r CCCC CCCC CCCC CCCC CCCC CCCC} \toprule
          & \rotatebox[origin=l]{90}{Portugal} & \rotatebox[origin=l]{90}{Croatia} & \rotatebox[origin=l]{90}{Austria} & \rotatebox[origin=l]{90}{Iceland} & \rotatebox[origin=l]{90}{Netherlands} & \rotatebox[origin=l]{90}{Spain} & \rotatebox[origin=l]{90}{Italy} & \rotatebox[origin=l]{90}{Czechia} & \rotatebox[origin=l]{90}{Poland} & \rotatebox[origin=l]{90}{Macedonia} & \rotatebox[origin=l]{90}{Qatar} & \rotatebox[origin=l]{90}{Brazil} & \rotatebox[origin=l]{90}{Argentina} & \rotatebox[origin=l]{90}{Cuba} & \rotatebox[origin=l]{90}{Japan} & \rotatebox[origin=l]{90}{Algeria} & \rotatebox[origin=l]{90}{Bahrain} & \rotatebox[origin=l]{90}{Tunisia} & \rotatebox[origin=l]{90}{Chile} & \rotatebox[origin=l]{90}{Kuwait} & \rotatebox[origin=l]{90}{Cape Verde} & \rotatebox[origin=l]{90}{Guinea} & \rotatebox[origin=l]{90}{United States} & \rotatebox[origin=l]{90}{Switzerland} \\ \bottomrule
    Denmark & \cellcolor{ForestGreen!4.70765714285714} 16.6 & X     & X     & \cellcolor{ForestGreen!4.82971428571429} 16.7 & \cellcolor{ForestGreen!4.77554285714286} 16.7 & \cellcolor{ForestGreen!4.74674285714286} 16.7 & \cellcolor{ForestGreen!4.69874285714286} 16.6 & \cellcolor{ForestGreen!4.81302857142857} 16.7 & \cellcolor{red!2.1064} 12.2 & \cellcolor{red!2.3296} 12.2 & \cellcolor{red!2.2192} 12.2 & \cellcolor{red!1.6808} 12.3 & \cellcolor{red!1.8192} 12.3 & \cellcolor{red!2.1624} 12.2 & \cellcolor{red!2.4472} 12.2 & \cellcolor{ForestGreen!2.10925714285714} 14.3 & \cellcolor{ForestGreen!1.10262857142857} 13.5 & \cellcolor{ForestGreen!1.96034285714286} 14.2 & \cellcolor{ForestGreen!1.11462857142857} 13.5 & \cellcolor{ForestGreen!0.933257142857143} 13.3 & \cellcolor{ForestGreen!2.02674285714286} 14.3 & \cellcolor{ForestGreen!2.05382857142857} 14.3 & \cellcolor{ForestGreen!1.26651428571429} 13.6 & \cellcolor{red!73.2056} 3.3 \\
    France & X     & \cellcolor{ForestGreen!42.7832} 49.9 & \cellcolor{ForestGreen!42.9310857142857} 50.1 & X     & X     & X     & X     & X     & \cellcolor{red!1.9144} 12.3 & \cellcolor{red!1.5208} 12.3 & \cellcolor{red!2.0184} 12.2 & \cellcolor{red!2.016} 12.2 & \cellcolor{red!2.0016} 12.2 & \cellcolor{red!2.3264} 12.2 & \cellcolor{red!2.2024} 12.2 & \cellcolor{ForestGreen!2} 14.3 & \cellcolor{ForestGreen!1.00628571428571} 13.4 & \cellcolor{ForestGreen!2.06011428571429} 14.3 & \cellcolor{ForestGreen!1.10285714285714} 13.5 & \cellcolor{ForestGreen!1.11508571428571} 13.5 & \cellcolor{ForestGreen!2.05725714285714} 14.3 & \cellcolor{ForestGreen!2.02022857142857} 14.3 & \cellcolor{ForestGreen!1.14742857142857} 13.5 & \cellcolor{red!73.5648} 3.3 \\
    Sweden & \cellcolor{ForestGreen!4.81942857142857} 16.7 & X     & X     & \cellcolor{ForestGreen!4.74057142857143} 16.6 & \cellcolor{ForestGreen!4.784} 16.7 & \cellcolor{ForestGreen!4.7416} 16.6 & \cellcolor{ForestGreen!4.7816} 16.7 & \cellcolor{ForestGreen!4.70422857142857} 16.6 & \cellcolor{red!1.616} 12.3 & \cellcolor{red!1.8336} 12.3 & \cellcolor{red!1.932} 12.3 & \cellcolor{red!2.0456} 12.2 & \cellcolor{red!1.9128} 12.3 & \cellcolor{red!2.0976} 12.2 & \cellcolor{red!2.456} 12.2 & \cellcolor{ForestGreen!1.9848} 14.2 & \cellcolor{ForestGreen!1.00342857142857} 13.4 & \cellcolor{ForestGreen!2.03257142857143} 14.3 & \cellcolor{ForestGreen!1.10182857142857} 13.5 & \cellcolor{ForestGreen!1.05565714285714} 13.4 & \cellcolor{ForestGreen!2.04582857142857} 14.3 & \cellcolor{ForestGreen!2.0544} 14.3 & \cellcolor{ForestGreen!1.18777142857143} 13.5 & \cellcolor{red!73.3704} 3.3 \\
    Germany & \cellcolor{ForestGreen!4.81611428571429} 16.7 & X     & X     & \cellcolor{ForestGreen!4.71554285714286} 16.6 & \cellcolor{ForestGreen!4.82331428571429} 16.7 & \cellcolor{ForestGreen!4.7152} 16.6 & \cellcolor{ForestGreen!4.77862857142857} 16.7 & \cellcolor{ForestGreen!4.72262857142857} 16.6 & \cellcolor{red!2.2952} 12.2 & \cellcolor{red!2.1328} 12.2 & \cellcolor{red!1.9832} 12.3 & \cellcolor{red!1.7944} 12.3 & \cellcolor{red!2.016} 12.2 & \cellcolor{red!2.092} 12.2 & \cellcolor{red!2.024} 12.2 & \cellcolor{ForestGreen!2.04822857142857} 14.3 & \cellcolor{ForestGreen!1.07005714285714} 13.4 & \cellcolor{ForestGreen!2.10571428571429} 14.3 & \cellcolor{ForestGreen!1.06022857142857} 13.4 & \cellcolor{ForestGreen!1.04034285714286} 13.4 & \cellcolor{ForestGreen!1.99942857142857} 14.2 & \cellcolor{ForestGreen!2.0008} 14.3 & \cellcolor{ForestGreen!1.18605714285714} 13.5 & \cellcolor{red!73.2384} 3.3 \\
    Hungary & \cellcolor{ForestGreen!4.7112} 16.6 & X     & X     & \cellcolor{ForestGreen!4.76971428571429} 16.7 & \cellcolor{ForestGreen!4.7408} 16.6 & \cellcolor{ForestGreen!4.80011428571429} 16.7 & \cellcolor{ForestGreen!4.75714285714286} 16.7 & \cellcolor{ForestGreen!4.79245714285714} 16.7 & \cellcolor{red!2.192} 12.2 & \cellcolor{red!2.1432} 12.2 & \cellcolor{red!1.7792} 12.3 & \cellcolor{red!1.8392} 12.3 & \cellcolor{red!2.2344} 12.2 & \cellcolor{red!2.46} 12.2 & \cellcolor{red!1.4328} 12.3 & \cellcolor{ForestGreen!2.01154285714286} 14.3 & \cellcolor{ForestGreen!1.06194285714286} 13.4 & \cellcolor{ForestGreen!2.00994285714286} 14.3 & \cellcolor{ForestGreen!1.02571428571429} 13.4 & \cellcolor{ForestGreen!1.02114285714286} 13.4 & \cellcolor{ForestGreen!2.0712} 14.3 & \cellcolor{ForestGreen!2.08925714285714} 14.3 & \cellcolor{ForestGreen!1.20914285714286} 13.6 & \cellcolor{red!73.4184} 3.3 \\
    Slovenia & \cellcolor{ForestGreen!4.80182857142857} 16.7 & X     & X     & \cellcolor{ForestGreen!4.7528} 16.7 & \cellcolor{ForestGreen!4.7024} 16.6 & \cellcolor{ForestGreen!4.76365714285714} 16.7 & \cellcolor{ForestGreen!4.76548571428571} 16.7 & \cellcolor{ForestGreen!4.78525714285714} 16.7 & \cellcolor{red!2.0632} 12.2 & \cellcolor{red!2.256} 12.2 & \cellcolor{red!2.2248} 12.2 & \cellcolor{red!2.2512} 12.2 & \cellcolor{red!2.2488} 12.2 & \cellcolor{red!1.4584} 12.3 & \cellcolor{red!2.2768} 12.2 & \cellcolor{ForestGreen!2.11131428571429} 14.3 & \cellcolor{ForestGreen!1.02708571428571} 13.4 & \cellcolor{ForestGreen!2.01737142857143} 14.3 & \cellcolor{ForestGreen!1.09805714285714} 13.5 & \cellcolor{ForestGreen!1.06274285714286} 13.4 & \cellcolor{ForestGreen!2.06971428571429} 14.3 & \cellcolor{ForestGreen!2.05702857142857} 14.3 & \cellcolor{ForestGreen!1.1656} 13.5 & \cellcolor{red!73.4832} 3.3 \\
    Norway & \cellcolor{ForestGreen!4.7152} 16.6 & X     & X     & \cellcolor{ForestGreen!4.76308571428571} 16.7 & \cellcolor{ForestGreen!4.74537142857143} 16.7 & \cellcolor{ForestGreen!4.80411428571429} 16.7 & \cellcolor{ForestGreen!4.78982857142857} 16.7 & \cellcolor{ForestGreen!4.75382857142857} 16.7 & \cellcolor{red!2.0408} 12.2 & \cellcolor{red!1.9488} 12.3 & \cellcolor{red!2.2888} 12.2 & \cellcolor{red!2.3048} 12.2 & \cellcolor{red!1.8824} 12.3 & \cellcolor{red!1.7816} 12.3 & \cellcolor{red!1.8968} 12.3 & \cellcolor{ForestGreen!2.02057142857143} 14.3 & \cellcolor{ForestGreen!1.04982857142857} 13.4 & \cellcolor{ForestGreen!2.09965714285714} 14.3 & \cellcolor{ForestGreen!1.09394285714286} 13.5 & \cellcolor{ForestGreen!1.04628571428571} 13.4 & \cellcolor{ForestGreen!2.01554285714286} 14.3 & \cellcolor{ForestGreen!2.01017142857143} 14.3 & \cellcolor{ForestGreen!1.14777142857143} 13.5 & \cellcolor{red!73.2424} 3.3 \\
    Egypt & X     & \cellcolor{ForestGreen!42.9310857142857} 50.1 & \cellcolor{ForestGreen!42.7832} 49.9 & X     & X     & X     & X     & X     & \cellcolor{ForestGreen!2.03257142857143} 14.3 & \cellcolor{ForestGreen!2.02354285714286} 14.3 & \cellcolor{ForestGreen!2.06365714285714} 14.3 & \cellcolor{ForestGreen!1.99028571428571} 14.2 & \cellcolor{ForestGreen!2.01645714285714} 14.3 & \cellcolor{ForestGreen!2.05405714285714} 14.3 & \cellcolor{ForestGreen!2.10514285714286} 14.3 & X     & \cellcolor{red!51.2488} 6.1 & X     & \cellcolor{red!53.1808} 5.9 & \cellcolor{red!50.9216} 6.1 & X     & X     & \cellcolor{red!58.172} 5.2 & \cellcolor{ForestGreen!73.3604571428571} 76.7 \\
    Portugal & X     & X     & X     & X     & X     & X     & X     & X     & \cellcolor{red!2.2024} 12.2 & \cellcolor{red!1.9192} 12.3 & \cellcolor{red!1.8312} 12.3 & \cellcolor{red!2.4792} 12.2 & \cellcolor{red!2.0672} 12.2 & \cellcolor{red!1.8152} 12.3 & \cellcolor{red!2.1952} 12.2 & \cellcolor{ForestGreen!2.0728} 14.3 & \cellcolor{ForestGreen!1.07348571428571} 13.4 & \cellcolor{ForestGreen!2.01405714285714} 14.3 & \cellcolor{ForestGreen!1.07622857142857} 13.4 & \cellcolor{ForestGreen!1.01348571428571} 13.4 & \cellcolor{ForestGreen!2.0528} 14.3 & \cellcolor{ForestGreen!2.03382857142857} 14.3 & \cellcolor{ForestGreen!1.2088} 13.6 & \cellcolor{red!73.3088} 3.3 \\
    Croatia & X     & X     & X     & X     & X     & X     & X     & X     & \cellcolor{ForestGreen!0.869828571428571} 13.3 & \cellcolor{ForestGreen!0.886971428571429} 13.3 & \cellcolor{ForestGreen!0.897942857142857} 13.3 & \cellcolor{ForestGreen!0.849142857142857} 13.2 & \cellcolor{ForestGreen!0.8504} 13.2 & \cellcolor{ForestGreen!0.893485714285714} 13.3 & \cellcolor{ForestGreen!0.919885714285714} 13.3 & \cellcolor{red!43.1736} 7.1 & \cellcolor{red!22.212} 9.7 & \cellcolor{red!42.8888} 7.1 & \cellcolor{red!22.8408} 9.6 & \cellcolor{red!21.5896} 9.8 & \cellcolor{red!43.1032} 7.1 & \cellcolor{red!42.7088} 7.2 & \cellcolor{red!25.172} 9.4 & \cellcolor{ForestGreen!31.5021714285714} 40.1 \\
    Austria & X     & X     & X     & X     & X     & X     & X     & X     & \cellcolor{ForestGreen!0.889257142857143} 13.3 & \cellcolor{ForestGreen!0.919314285714286} 13.3 & \cellcolor{ForestGreen!0.877371428571429} 13.3 & \cellcolor{ForestGreen!0.853142857142857} 13.2 & \cellcolor{ForestGreen!0.880114285714286} 13.3 & \cellcolor{ForestGreen!0.828228571428571} 13.2 & \cellcolor{ForestGreen!0.870628571428571} 13.3 & \cellcolor{red!42.8264} 7.1 & \cellcolor{red!21.9928} 9.8 & \cellcolor{red!42.6904} 7.2 & \cellcolor{red!22.62} 9.7 & \cellcolor{red!21.5264} 9.8 & \cellcolor{red!42.496} 7.2 & \cellcolor{red!43.1496} 7.1 & \cellcolor{red!24.968} 9.4 & \cellcolor{ForestGreen!31.3490285714286} 39.9 \\
    Iceland & X     & X     & X     & X     & X     & X     & X     & X     & \cellcolor{red!2.096} 12.2 & \cellcolor{red!2.5032} 12.2 & \cellcolor{red!2.3296} 12.2 & \cellcolor{red!1.8688} 12.3 & \cellcolor{red!1.864} 12.3 & \cellcolor{red!1.8624} 12.3 & \cellcolor{red!1.7848} 12.3 & \cellcolor{ForestGreen!2.04411428571429} 14.3 & \cellcolor{ForestGreen!1.03268571428571} 13.4 & \cellcolor{ForestGreen!2.02537142857143} 14.3 & \cellcolor{ForestGreen!1.06765714285714} 13.4 & \cellcolor{ForestGreen!0.991314285714286} 13.4 & \cellcolor{ForestGreen!2.04057142857143} 14.3 & \cellcolor{ForestGreen!2.06445714285714} 14.3 & \cellcolor{ForestGreen!1.22148571428571} 13.6 & \cellcolor{red!73.1048} 3.4 \\
    Netherlands & X     & X     & X     & X     & X     & X     & X     & X     & \cellcolor{red!1.6504} 12.3 & \cellcolor{red!2.156} 12.2 & \cellcolor{red!2.448} 12.2 & \cellcolor{red!2} 12.3 & \cellcolor{red!1.5488} 12.3 & \cellcolor{red!2.656} 12.2 & \cellcolor{red!1.988} 12.3 & \cellcolor{ForestGreen!2.06388571428571} 14.3 & \cellcolor{ForestGreen!1.07942857142857} 13.4 & \cellcolor{ForestGreen!2.03702857142857} 14.3 & \cellcolor{ForestGreen!1.10571428571429} 13.5 & \cellcolor{ForestGreen!1.048} 13.4 & \cellcolor{ForestGreen!1.99051428571429} 14.2 & \cellcolor{ForestGreen!1.984} 14.2 & \cellcolor{ForestGreen!1.24102857142857} 13.6 & \cellcolor{red!73.4} 3.3 \\
    Spain & X     & X     & X     & X     & X     & X     & X     & X     & \cellcolor{red!1.8432} 12.3 & \cellcolor{red!1.8856} 12.3 & \cellcolor{red!2.016} 12.2 & \cellcolor{red!2.388} 12.2 & \cellcolor{red!1.6608} 12.3 & \cellcolor{red!2.216} 12.2 & \cellcolor{red!2.3792} 12.2 & \cellcolor{ForestGreen!2.05554285714286} 14.3 & \cellcolor{ForestGreen!1.03348571428571} 13.4 & \cellcolor{ForestGreen!2.13291428571429} 14.4 & \cellcolor{ForestGreen!1.03554285714286} 13.4 & \cellcolor{ForestGreen!1.05245714285714} 13.4 & \cellcolor{ForestGreen!1.9968} 14.2 & \cellcolor{ForestGreen!2.07794285714286} 14.3 & \cellcolor{ForestGreen!1.17211428571429} 13.5 & \cellcolor{red!73.5088} 3.3 \\
    Italy & X     & X     & X     & X     & X     & X     & X     & X     & \cellcolor{red!2.4568} 12.2 & \cellcolor{red!2.22} 12.2 & \cellcolor{red!1.6768} 12.3 & \cellcolor{red!1.2376} 12.3 & \cellcolor{red!2.536} 12.2 & \cellcolor{red!1.8088} 12.3 & \cellcolor{red!2.0952} 12.2 & \cellcolor{ForestGreen!2.00445714285714} 14.3 & \cellcolor{ForestGreen!1.0648} 13.4 & \cellcolor{ForestGreen!2.0312} 14.3 & \cellcolor{ForestGreen!1.08} 13.4 & \cellcolor{ForestGreen!1.03371428571429} 13.4 & \cellcolor{ForestGreen!2.04262857142857} 14.3 & \cellcolor{ForestGreen!2.05462857142857} 14.3 & \cellcolor{ForestGreen!1.16377142857143} 13.5 & \cellcolor{red!73.2952} 3.3 \\
    Czechia & X     & X     & X     & X     & X     & X     & X     & X     & \cellcolor{red!2.0648} 12.2 & \cellcolor{red!1.96} 12.3 & \cellcolor{red!2.1256} 12.2 & \cellcolor{red!1.9424} 12.3 & \cellcolor{red!2.4368} 12.2 & \cellcolor{red!1.6936} 12.3 & \cellcolor{red!2.0912} 12.2 & \cellcolor{ForestGreen!2.04491428571429} 14.3 & \cellcolor{ForestGreen!1.03108571428571} 13.4 & \cellcolor{ForestGreen!1.98502857142857} 14.2 & \cellcolor{ForestGreen!1.12925714285714} 13.5 & \cellcolor{ForestGreen!1.02045714285714} 13.4 & \cellcolor{ForestGreen!2.10514285714286} 14.3 & \cellcolor{ForestGreen!2.05062857142857} 14.3 & \cellcolor{ForestGreen!1.15565714285714} 13.5 & \cellcolor{red!73.3408} 3.3 \\
    Poland & \cellcolor{red!2.2024} 12.2 & \cellcolor{ForestGreen!0.869828571428571} 13.3 & \cellcolor{ForestGreen!0.889257142857143} 13.3 & \cellcolor{red!2.096} 12.2 & \cellcolor{red!1.6504} 12.3 & \cellcolor{red!1.8432} 12.3 & \cellcolor{red!2.4568} 12.2 & \cellcolor{red!2.0648} 12.2 & X     & X     & X     & X     & X     & X     & X     & X     & \cellcolor{ForestGreen!4.2416} 16.2 & \cellcolor{red!2.9992} 12.1 & \cellcolor{ForestGreen!3.24251428571429} 15.3 & \cellcolor{ForestGreen!4.28628571428571} 16.3 & \cellcolor{red!2.76} 12.2 & \cellcolor{red!2.2208} 12.2 & \cellcolor{ForestGreen!0.663428571428571} 13.1 & \cellcolor{red!79.0568} 2.6 \\
    Macedonia & \cellcolor{red!1.9192} 12.3 & \cellcolor{ForestGreen!0.886971428571429} 13.3 & \cellcolor{ForestGreen!0.919314285714286} 13.3 & \cellcolor{red!2.5032} 12.2 & \cellcolor{red!2.156} 12.2 & \cellcolor{red!1.8856} 12.3 & \cellcolor{red!2.22} 12.2 & \cellcolor{red!1.96} 12.3 & X     & X     & X     & X     & X     & X     & X     & X     & \cellcolor{ForestGreen!4.26274285714286} 16.2 & \cellcolor{red!2.3752} 12.2 & \cellcolor{ForestGreen!3.20445714285714} 15.3 & \cellcolor{ForestGreen!4.25257142857143} 16.2 & \cellcolor{red!2.3856} 12.2 & \cellcolor{red!2.6256} 12.2 & \cellcolor{ForestGreen!0.625028571428571} 13 & \cellcolor{red!79.0272} 2.6 \\
    Qatar & \cellcolor{red!1.8312} 12.3 & \cellcolor{ForestGreen!0.897942857142857} 13.3 & \cellcolor{ForestGreen!0.877371428571429} 13.3 & \cellcolor{red!2.3296} 12.2 & \cellcolor{red!2.448} 12.2 & \cellcolor{red!2.016} 12.2 & \cellcolor{red!1.6768} 12.3 & \cellcolor{red!2.1256} 12.2 & X     & X     & X     & X     & X     & X     & X     & X     & X     & \cellcolor{ForestGreen!5.08365714285714} 16.9 & \cellcolor{ForestGreen!4.30217142857143} 16.3 & X     & \cellcolor{ForestGreen!5.12765714285714} 17 & \cellcolor{ForestGreen!5.12708571428571} 17 & \cellcolor{ForestGreen!2.28617142857143} 14.5 & \cellcolor{ForestGreen!6.64468571428571} 18.3 \\
    Brazil & \cellcolor{red!2.4792} 12.2 & \cellcolor{ForestGreen!0.849142857142857} 13.2 & \cellcolor{ForestGreen!0.853142857142857} 13.2 & \cellcolor{red!1.8688} 12.3 & \cellcolor{red!2} 12.3 & \cellcolor{red!2.388} 12.2 & \cellcolor{red!1.2376} 12.3 & \cellcolor{red!1.9424} 12.3 & X     & X     & X     & X     & X     & X     & X     & X     & \cellcolor{ForestGreen!2.00731428571429} 14.3 & \cellcolor{ForestGreen!1.76788571428571} 14 & X     & \cellcolor{ForestGreen!2.00902857142857} 14.3 & \cellcolor{ForestGreen!1.73371428571429} 14 & \cellcolor{ForestGreen!1.68594285714286} 14 & \cellcolor{red!5.188} 11.9 & \cellcolor{ForestGreen!5.82297142857143} 17.6 \\
    Argentina & \cellcolor{red!2.0672} 12.2 & \cellcolor{ForestGreen!0.8504} 13.2 & \cellcolor{ForestGreen!0.880114285714286} 13.3 & \cellcolor{red!1.864} 12.3 & \cellcolor{red!1.5488} 12.3 & \cellcolor{red!1.6608} 12.3 & \cellcolor{red!2.536} 12.2 & \cellcolor{red!2.4368} 12.2 & X     & X     & X     & X     & X     & X     & X     & X     & \cellcolor{ForestGreen!2.05382857142857} 14.3 & \cellcolor{ForestGreen!1.78457142857143} 14.1 & X     & \cellcolor{ForestGreen!1.97382857142857} 14.2 & \cellcolor{ForestGreen!1.72068571428571} 14 & \cellcolor{ForestGreen!1.65074285714286} 13.9 & \cellcolor{red!5.48} 11.8 & \cellcolor{ForestGreen!5.88491428571429} 17.6 \\
    Cuba  & \cellcolor{red!1.8152} 12.3 & \cellcolor{ForestGreen!0.893485714285714} 13.3 & \cellcolor{ForestGreen!0.828228571428571} 13.2 & \cellcolor{red!1.8624} 12.3 & \cellcolor{red!2.656} 12.2 & \cellcolor{red!2.216} 12.2 & \cellcolor{red!1.8088} 12.3 & \cellcolor{red!1.6936} 12.3 & X     & X     & X     & X     & X     & X     & X     & X     & \cellcolor{ForestGreen!1.772} 14.1 & \cellcolor{ForestGreen!1.33097142857143} 13.7 & \cellcolor{ForestGreen!0.946057142857143} 13.3 & \cellcolor{ForestGreen!1.728} 14 & \cellcolor{ForestGreen!1.32388571428571} 13.7 & \cellcolor{ForestGreen!1.38914285714286} 13.7 & X     & \cellcolor{ForestGreen!5.79565714285714} 17.6 \\
    Japan & \cellcolor{red!2.1952} 12.2 & \cellcolor{ForestGreen!0.919885714285714} 13.3 & \cellcolor{ForestGreen!0.870628571428571} 13.3 & \cellcolor{red!1.7848} 12.3 & \cellcolor{red!1.988} 12.3 & \cellcolor{red!2.3792} 12.2 & \cellcolor{red!2.0952} 12.2 & \cellcolor{red!2.0912} 12.2 & X     & X     & X     & X     & X     & X     & X     & X     & X     & \cellcolor{ForestGreen!5.0864} 17 & \cellcolor{ForestGreen!4.26937142857143} 16.2 & X     & \cellcolor{ForestGreen!5.11485714285714} 17 & \cellcolor{ForestGreen!5.12514285714286} 17 & \cellcolor{ForestGreen!2.25302857142857} 14.5 & \cellcolor{ForestGreen!6.72262857142857} 18.4 \\
    Algeria & \cellcolor{ForestGreen!2.0728} 14.3 & \cellcolor{red!43.1736} 7.1 & \cellcolor{red!42.8264} 7.1 & \cellcolor{ForestGreen!2.04411428571429} 14.3 & \cellcolor{ForestGreen!2.06388571428571} 14.3 & \cellcolor{ForestGreen!2.05554285714286} 14.3 & \cellcolor{ForestGreen!2.00445714285714} 14.3 & \cellcolor{ForestGreen!2.04491428571429} 14.3 & X     & X     & X     & X     & X     & X     & X     & X     & \cellcolor{ForestGreen!14.2339428571429} 25 & X     & \cellcolor{ForestGreen!12.6068571428571} 23.5 & \cellcolor{ForestGreen!14.3217142857143} 25 & X     & X     & \cellcolor{ForestGreen!9.98205714285714} 21.2 & \cellcolor{red!58.012} 5.2 \\ \toprule
    \end{tabularx}
\end{tiny}
    \begin{tablenotes} \footnotesize
\item
Abbreviations: Czechia = Czech Republic; Macedonia = North Macedonia.
\item
X represents a pair of teams that cannot play in the preliminary group stage.
\item
The numbers show percentages rounded to one decimal place.
\item
\textcolor{ForestGreen}{Green} (\textcolor{red}{Red}) colour means that the draw procedure implies a higher (lower) probability than the equal probability of 1/8 = 12.5\%.
\item
Darker colour indicates a higher distortion.
    \end{tablenotes}
\end{threeparttable}
}
\end{subtable}
\end{table}

%\end{document}

%% file: Group_draw_attractiveness_equal_treatment_handball.bbl
\begin{thebibliography}{}

\bibitem[Boczo{\'n} and Wilson, 2023]{BoczonWilson2023}
Boczo{\'n}, M. and Wilson, A.~J. (2023).
\newblock Goals, constraints, and transparently fair assignments: A field study
  of randomization design in the {UEFA} {C}hampions {L}eague.
\newblock {\em Management Science}, 69(6):3474--3491.

\bibitem[Cea et~al., 2020]{CeaDuranGuajardoSureSiebertZamorano2020}
Cea, S., Dur{\'a}n, G., Guajardo, M., Saur{\'e}, D., Siebert, J., and Zamorano,
  G. (2020).
\newblock An analytics approach to the {FIFA} ranking procedure and the {W}orld
  {C}up final draw.
\newblock {\em Annals of Operations Research}, 286(1-2):119--146.

\bibitem[Csat{\'o}, 2021]{Csato2021a}
Csat{\'o}, L. (2021).
\newblock {\em Tournament Design: How Operations Research Can Improve Sports
  Rules}.
\newblock Palgrave Pivots in Sports Economics. Palgrave Macmillan, Cham,
  Switzerland.

\bibitem[Csat{\'o}, 2022]{Csato2022a}
Csat{\'o}, L. (2022).
\newblock Quantifying incentive (in)compatibility: {A} case study from sports.
\newblock {\em European Journal of Operational Research}, 302(2):717--726.

\bibitem[Csat{\'o}, 2023a]{Csato2023d}
Csat{\'o}, L. (2023a).
\newblock Group draw with unknown qualified teams: {A} lesson from the 2022
  {FIFA} {W}orld {C}up.
\newblock {\em International Journal of Sports Science \& Coaching},
  18(2):539--551.

\bibitem[Csat{\'o}, 2023b]{Csato2023c}
Csat{\'o}, L. (2023b).
\newblock Quantifying the unfairness of the 2018 {FIFA} {W}orld {C}up
  qualification.
\newblock {\em International Journal of Sports Science \& Coaching},
  18(1):183--196.

\bibitem[Csat{\'o}, 2024a]{Csato2024i}
Csat{\'o}, L. (2024a).
\newblock The fairness of transparent constrained assignments: The role of draw
  order.
\newblock Manuscript. {DOI}:
  \href{https://doi.org/10.48550/arXiv.2109.13785}{10.48550/arXiv.2109.13785}.

\bibitem[Csat{\'o}, 2024b]{Csato2024h}
Csat{\'o}, L. (2024b).
\newblock Random matching in balanced bipartite graphs: The (un)fairness of
  draw mechanisms used in sports.
\newblock Manuscript. {DOI}:
  \href{https://doi.org/10.48550/arXiv.2303.09274}{10.48550/arXiv.2303.09274}.

\bibitem[Csat{\'o}, 2025]{Csato2025c}
Csat{\'o}, L. (2025).
\newblock The fairness of the group draw for the {FIFA} {W}orld {C}up.
\newblock {\em International Journal of Sports Science \& Coaching}, in press.
\newblock {DOI}:
  \href{https://doi.org/10.1177/17479541241300219}{10.1177/17479541241300219}.

\bibitem[Csat{\'o} et~al., 2025]{CsatoKissSzadoczki2025}
Csat{\'o}, L., Kiss, L.~M., and Sz{\'a}doczki, {\relax Zs}. (2025).
\newblock The allocation of {FIFA} {W}orld {C}up slots based on the ranking of
  confederations.
\newblock {\em Annals of Operations Research}, 344(1):153--173.

\bibitem[Devriesere et~al., 2024]{DevriesereCsatoGoossens2024}
Devriesere, K., Csat{\'o}, L., and Goossens, D. (2024).
\newblock Tournament design: A review from an operational research perspective.
\newblock {\em European Journal of Operational Research}, in press.
\newblock {DOI}:
  \href{https://doi.org/10.1016/j.ejor.2024.10.044}{10.1016/j.ejor.2024.10.044}.

\bibitem[FIBA, 2019]{FIBA2019}
FIBA (2019).
\newblock Procedure for {FIBA} {B}asketball {W}orld {C}up 2019 {D}raw.
\newblock 15 March.
  \url{https://www.fiba.basketball/basketballworldcup/2019/news/procedure-for-fiba-basketball-world-cup-2019-draw}.

\bibitem[FIBA, 2023]{FIBA2023}
FIBA (2023).
\newblock The {FIBA} {B}asketball {W}orld {C}up 2023 draw principles explained.
\newblock 21 April.
  \url{https://www.fiba.basketball/basketballworldcup/2023/news/the-fiba-basketball-world-cup-2023-draw-principles-explained}.

\bibitem[FIFA, 2017]{FIFA2017c}
FIFA (2017).
\newblock Close-up on {F}inal {D}raw procedures.
\newblock 27 November.
  \url{https://web.archive.org/web/20171127150059/http://www.fifa.com/worldcup/news/y=2017/m=11/news=close-up-on-final-draw-procedures-2921440.html}.

\bibitem[FIFA, 2022]{FIFA2022a}
FIFA (2022).
\newblock {\em {D}raw procedures. {FIFA} {W}orld {C}up {Q}atar
  2022\textsuperscript{{TM}}}.
\newblock
  \url{https://digitalhub.fifa.com/m/2ef762dcf5f577c6/original/Portrait-Master-Template.pdf}.

\bibitem[Goossens et~al., 2020]{GoossensYiVanBulck2020}
Goossens, D., Yi, X., and Van~Bulck, D. (2020).
\newblock Fairness trade-offs in sports timetabling.
\newblock In Ley, C. and Dominicy, Y., editors, {\em Science Meets Sports: When
  Statistics Are More Than Numbers}, pages 213--244. Cambridge Scholars
  Publishing, Newcastle upon Tyne, United Kingdom.

\bibitem[Guyon, 2014]{Guyon2014a}
Guyon, J. (2014).
\newblock Rethinking the {FIFA} {W}orld {C}up\textsuperscript{{TM}} final draw.
\newblock Manuscript. {DOI}:
  \href{http://dx.doi.org/10.2139/ssrn.2424376}{10.2139/ssrn.2424376}.

\bibitem[Guyon, 2015]{Guyon2015a}
Guyon, J. (2015).
\newblock Rethinking the {FIFA} {W}orld {C}up\textsuperscript{{TM}} final draw.
\newblock {\em Journal of Quantitative Analysis in Sports}, 11(3):169--182.

\bibitem[Guyon, 2018]{Guyon2018d}
Guyon, J. (2018).
\newblock Pourquoi la {C}oupe du monde est plus \'equitable cette ann\'ee.
\newblock {\em The Conversation}.
\newblock 13 June.
  \url{https://theconversation.com/pourquoi-la-coupe-du-monde-est-plus-equitable-cette-annee-97948}.

\bibitem[Guyon and Meunier, 2023]{GuyonMeunier2023}
Guyon, J. and Meunier, F. (2023).
\newblock {UEFA} draws: probability calculator.
\newblock 13 December.
  \url{https://julienguyon.github.io/UEFA-draws/ExactProbabilities_ChampionsLeague_RoundOf16_text_website-1.pdf}.

\bibitem[Guyon et~al., 2024]{GuyonMeunierBenSalemBuchholtzerTanre2024}
Guyon, J., Meunier, F., Ben~Salem, A., Buchholtzer, T., and Tanr\'e, M. (2024).
\newblock The drawing of the teams in the future {C}hampions {L}eague.
\newblock Academic report.
  \url{https://julienguyon.pythonanywhere.com/uefa_app/}.

\bibitem[Gyimesi, 2024]{Gyimesi2024}
Gyimesi, A. (2024).
\newblock Competitive balance in the post-2024 {C}hampions {L}eague and the
  {E}uropean {S}uper {L}eague: {A} simulation study.
\newblock {\em Journal of Sports Economics}, 25(6):707--734.

\bibitem[Hall, 1935]{Hall1935}
Hall, P. (1935).
\newblock On representatives of subsets.
\newblock {\em Journal of the London Mathematical Society}, 1(1):26--30.

\bibitem[Herfindahl, 1950]{Herfindahl1950}
Herfindahl, O.~C. (1950).
\newblock {\em Concentration in the Steel Industry}.
\newblock PhD thesis, Columbia University, New York.

\bibitem[Hirschman, 1945]{Hirschman1945}
Hirschman, A.~O. (1945).
\newblock {\em National {P}ower and the {S}tructure of {F}oreign {T}rade}.
\newblock University of California Press, Berkeley and Los Angeles, California,
  USA.

\bibitem[IHF, 2025]{IHF2025}
IHF (2025).
\newblock Pots revealed for the 2025 {IHF} {M}en's {W}orld {C}hampionship draw.
\newblock 24 May.
  \url{https://www.ihf.info/media-center/news/pots-revealed-2025-ihf-mens-world-championship-draw}.

\bibitem[Kiesl, 2013]{Kiesl2013}
Kiesl, H. (2013).
\newblock Match me if you can. {M}athematische {G}edanken zur
  {C}hampions-{L}eague-{A}chtelfinalauslosung.
\newblock {\em Mitteilungen der Deutschen Mathematiker-Vereinigung},
  21(2):84--88.

\bibitem[Kl{\"o}{\ss}ner and Becker, 2013]{KlossnerBecker2013}
Kl{\"o}{\ss}ner, S. and Becker, M. (2013).
\newblock Odd odds: The {UEFA} {C}hampions {L}eague {R}ound of 16 draw.
\newblock {\em Journal of Quantitative Analysis in Sports}, 9(3):249--270.

\bibitem[Krumer and Moreno-Ternero, 2023]{KrumerMoreno-Ternero2023}
Krumer, A. and Moreno-Ternero, J. (2023).
\newblock The allocation of additional slots for the {FIFA} {W}orld {C}up.
\newblock {\em Journal of Sports Economics}, 24(7):831--850.

\bibitem[Laliena and L{\'o}pez, 2019]{LalienaLopez2019}
Laliena, P. and L{\'o}pez, F.~J. (2019).
\newblock Fair draws for group rounds in sport tournaments.
\newblock {\em International Transactions in Operational Research},
  26(2):439--457.

\bibitem[Laliena and L{\'o}pez, 2025]{LalienaLopez2025}
Laliena, P. and L{\'o}pez, F.~J. (2025).
\newblock Draw procedures for balanced 3-team group rounds in sports
  competitions.
\newblock {\em Annals of Operations Research}, in press.
\newblock {DOI}:
  \href{https://doi.org/10.1007/s10479-025-06497-9}{10.1007/s10479-025-06497-9}.

\bibitem[Lasek and Gagolewski, 2018]{LasekGagolewski2018}
Lasek, J. and Gagolewski, M. (2018).
\newblock The efficacy of league formats in ranking teams.
\newblock {\em Statistical Modelling}, 18(5-6):411--435.

\bibitem[Owen and Owen, 2022]{OwenOwen2022}
Owen, P.~D. and Owen, C.~A. (2022).
\newblock Simulation evidence on {H}erfindahl-{H}irschman measures of
  competitive balance in professional sports leagues.
\newblock {\em Journal of the Operational Research Society}, 73(2):285--300.

\bibitem[Owen et~al., 2007]{OwenRyanWeatherston2007}
Owen, P.~D., Ryan, M., and Weatherston, C.~R. (2007).
\newblock Measuring competitive balance in professional team sports using the
  {H}erfindahl-{H}irschman index.
\newblock {\em Review of Industrial Organization}, 31(4):289--302.

\bibitem[Roberts and Rosenthal, 2024]{RobertsRosenthal2024}
Roberts, G.~O. and Rosenthal, J.~S. (2024).
\newblock Football group draw probabilities and corrections.
\newblock {\em Canadian Journal of Statistics}, 52(3):659--677.

\bibitem[Sauer et~al., 2024]{SauerCsehLenzner2024}
Sauer, P., Cseh, {\'A}., and Lenzner, P. (2024).
\newblock Improving ranking quality and fairness in {S}wiss-system chess
  tournaments.
\newblock {\em Journal of Quantitative Analysis in Sports}, 20(2):127--146.

\bibitem[Sziklai et~al., 2022]{SziklaiBiroCsato2022}
Sziklai, B.~R., Bir\'o, P., and Csat{\'o}, L. (2022).
\newblock The efficacy of tournament designs.
\newblock {\em Computers \& Operations Research}, 144:105821.

\bibitem[Tijms, 2015]{Tijms2015}
Tijms, H. (2015).
\newblock Teaching note -- {W}as the {C}hampions {L}eague draw rigged?
\newblock {\em Teaching Statistics}, 37(3):104--105.

\bibitem[van Ours, 2024]{vanOurs2024}
van Ours, J.~C. (2024).
\newblock Nontransitive patterns in long-term football rivalries.
\newblock {\em Journal of Sports Economics}, 25(7):802--826.

\bibitem[{Volleyball World}, 2024]{Volleyball2024}
{Volleyball World} (2024).
\newblock 2025 {M}en’s {W}orld {C}hampionship pools to be drawn on
  {S}eptember 14.
\newblock 4 September.
  \url{https://en.volleyballworld.com/volleyball/competitions/men-world-championship/news/2025-men-s-world-championship-pools-to-be-drawn-on-september-14}.

\bibitem[Wallace and Haigh, 2013]{WallaceHaigh2013}
Wallace, M. and Haigh, J. (2013).
\newblock Football and marriage -- and the {UEFA} draw.
\newblock {\em Significance}, 10(2):47--48.

\end{thebibliography}
